\documentclass{article}
\usepackage{amsmath}
\usepackage{amsfonts,amssymb}
\usepackage[dvips]{graphics}
\usepackage{graphicx}

\input amssym.def

\numberwithin{equation}{section}

\newcommand{\cHom}{\mathcal Hom\,}

\newcommand{\Vol}{{\rm V o l }}
\newcommand{\MC}{{\rm M C }}

\newcommand{\Hom}{{\rm H o m }}
\newcommand{\Lie}{{\rm L i e }}

\newcommand{\Der}{{\rm D e r }}

\newcommand{\tcL}{\widetilde{\cal L}}

\newcommand{\Omb}{\Om^{-\bul}}

\newcommand{\sCbu}{C^{\bullet+1}}
\newcommand{\Cbu}{C^{\bullet}}
\newcommand{\Cbd}{C_{\bullet}}

\newcommand{\Linf}{L_{\infty}}

\newcommand{\tpi}{\widetilde{\pi}}

\newcommand{\hA}{{A_{\hbar}}}

\newcommand{\tV}{\widetilde{V}}

\newcommand{\al}{{\alpha}}

\newcommand{\h}{{\hbar}}
\newcommand{\bul}{{\bullet}}

\newcommand{\mG}{{\mathfrak{G}}}

\newcommand{\Om}{{\Omega}}

\newcommand{\si}{{\sigma}}
\newcommand{\ga}{{\gamma}}

\newcommand{\ve}{{\varepsilon}}

\newcommand{\pa}{{\partial}}

\newcommand{\bs}{{\bf s}}

\newcommand{\cK}{{\cal K}}

\newcommand{\cL}{{\cal L}}
\newcommand{\cR}{{\cal R}}

\newcommand{\cG}{{\cal G}}

\newcommand{\cW}{{\cal W}}

\newcommand{\cO}{{\cal O}}

\newcommand{\bbK}{{\mathbb K}}

\newcommand{\bbR}{{\mathbb R}}

\newcommand{\te}{\theta}

\newcommand{\de}{{\delta}}
\newcommand{\D}{{\Delta}}

\newcommand{\tQ}{{\widetilde{Q}}}

\date{}
\newtheorem{defi}{Definition}

\newtheorem{teo}{Theorem}
\newtheorem{cor}{Corollary}

\newtheorem{pred}{Proposition}

\title{Formality theorems for Hochschild complexes
and their applications}

\author{V.A. Dolgushev, D.E. Tamarkin, and B.L. Tsygan}
\date{}

\begin{document}

\large

\maketitle

\begin{center}
{\it To Giovanni Felder on the occasion of his 50th birthday.}
\end{center}

\begin{abstract}
We give a popular introduction to formality theorems for
Hochschild complexes and their applications. We review
some of the recent results and prove that the truncated
Hochschild cochain complex of a polynomial
algebra is non-formal.
\end{abstract}

\tableofcontents

\section{Introduction}
The notion of formality was suggested in classical
article \cite{DGMS} of  P. Deligne, P. Griffiths, J. Morgan,
and D. Sullivan. In this article it was shown that the
de Rham algebra of a compact K\"ahler manifold $M$ is quasi-isomorphic
to the cohomology ring of $M$\,. Using the terminology suggested
in this article we can say that de Rham algebra of a compact
K\"ahler manifold $M$ is {\it formal} as a commutative algebra.

Around 1993-94 M. Kontsevich conjectured
(see \cite{K1}, \cite{Sasha-V}) that the
Hochschild cochain complex for the algebra of functions
on a smooth manifold is formal as a Lie algebra with
the Gerstenhaber bracket \cite{G}.
Then in 1997 M. Kontsevich proved \cite{K} this formality
conjecture for an arbitrary smooth manifold.
In the same paper he showed how this result
solves a long standing problem on the deformation quantization
\cite{Bayen}, \cite{Ber} of a Poisson manifold.

In 1998 the second author proposed a completely different
proof of Kontsevich's formality theorem \cite{Dima-Proof} for
the case of the affine space over an arbitrary field of
characteristic zero.
This approach is based on deep results such as the proof of Deligne's
conjecture on the Hochschild complex \cite{BF}, \cite{K-Soi}, \cite{M-Smith},
\cite{Dima-DG}, \cite{Sasha-Proof} and the formality theorem \cite{K-oper},
\cite{LV},  \cite{Dima-Disk} for the operad of little discs.

In 1999 A.S. Cattaneo and G. Felder described \cite{CF-path} how
Kontsevich's star-product formula as well as his formality theorem can
be obtained using the correlators of the Poisson sigma model
\cite{Ikeda}, \cite{Thomas}.

After M. Kontsevich's celebrated result \cite{K}
lots of interesting generalizations and applications of
the formality theorem for Hochschild cochain complex were
proposed. At this moment all these results can be put
under an umbrella of an independent mathematical topic.
In this paper we give a popular introduction to this
fascinating topic. We review some of the recent results
and give an example of a non-formal differential graded (DG) Lie algebra.
We hope that our introductory part is accessible to
graduate students who are interested in this topic.

The organization of the paper is as follows.
In the next section we illustrate the general concept
of formality with the example of a DG
Lie algebra. We recall Maurer-Cartan elements, the Goldman-Millson
groupoid and twisting procedure. We also discuss some consequences
of formality for a DG Lie algebra.
In Section \ref{alg-str}
we recall basic algebraic structures
on the Hochschild complexes of an associative algebra $A$\,.
In this section we also recall the Van den Bergh duality theorem
\cite{VB}.
Section \ref{F} begins with the formulation of
Kontsevich's formality theorem \cite{K} and its immediate
corollaries. Next we review the alternative approach \cite{Dima-Proof}
of the second author. Then we discuss formality theorems for
Hochschild and cyclic chains and the formality of the 
$\infty$-calculus algebra of Hochschild complexes.
We conclude Section \ref{F} with the formality theorems
for Hochschild complexes in the algebraic geometry setting.
In Section \ref{more-applic} we give a brief outline of few recent
applications of formality theorems
for Hochschild complexes.  Finally, in the concluding section
we give an example of a non-formal DG Lie algebra.

~\\
{\bf Notation.}
"DG" stands for differential graded.
$S_n$ denotes the symmetric group on $n$ letters.
By suspension $\bs V$ of a graded vector space
(or a cochain complex) $V$ we mean $\ve \otimes V$, where $\ve$ is a
one-dimensional vector space placed in degree
$+1$\,.  For a vector $v\in V$ we denote by $|v|$
its degree.
We use the Koszul rule of signs which says that
a transposition of any two homogeneous vectors $v_1$ and
$v_2$ yields the sign
$$
(-1)^{|v_1| \, |v_2|}\,.
$$
For a groupoid $\cG$ we denote by $\pi_0(\cG)$ the
set of isomorphism classes of its objects. This is exactly
the set of connected components of the classifying space
$B \cG$ of $\cG$\,.
We assume that the underlying field $\bbK$
has characteristic zero. $\h$ denotes the formal deformation
parameter.

~\\
{\bf Acknowledgment.}
We would like to thank J. Stasheff for discussions and for 
his useful comments on the first version of the manuscript. 
D.T. and B.T. are supported by
NSF grants. The work of V.D. is partially supported by the Grant
for Support of Scientific Schools NSh-3036.2008.2

\section{Formal versus non-formal}
\label{f-versus-nonf}
\subsection{General definitions}
Let $\cL$ be a DG Lie algebra over the field $\bbK$\,.
We say that a morphism $\mu : \cL \to \tcL $ is
a {\it quasi-isomorphism} if $\mu$ induces an isomorphism
on cohomology groups. In this case we use the tilde
over the arrow
$$
\mu : \cL \stackrel{\sim}{\to} \tcL\,.
$$

We call two DG Lie algebras
$\cL$ and $\tcL$ {\it quasi-isomorphic} if they can be connected
by a sequence of quasi-isomorphisms
$$
\cL\,  \stackrel{\sim}{\rightarrow} \,\bul\,
\stackrel{\sim}{\leftarrow} \,\bul\, \stackrel{\sim}{\rightarrow}
\,\bul\,
\dots \,\bul\, \stackrel{\sim}{\leftarrow} \, \tcL \,.
$$

For every DG Lie algebra $\cL$ its cohomology $H^{\bul}(\cL)$
is naturally a graded Lie algebra. We think of $H^{\bul}(\cL)$
as the DG Lie algebra with the zero differential.
\begin{defi}
\label{formality}
A DG Lie algebra $\cL$ is called {\rm formal} if
it is quasi-isomorphic to its cohomology
$H^{\bul}(\cL)$\,.
\end{defi}
{\bf Remark.} Similarly, we can talk
about formal or non-formal DG associative algebras, DG commutative
algebras and DG algebras of other types.

We would like to mention that, since our DG Lie
algebras are complexes of vector spaces over a field, we always
have a map of cochain complexes
$$
F_1: H^{\bul}(\cL) \to \cL
$$
which induces an isomorphism on the level of cohomology.
In other words, it is always possible to choose
a representative for each cohomology class in such
a way that this choice respects linearity.

It is obvious that, in general, $F_1$
is not compatible with the Lie brackets.
However, there exists a bilinear map
$$
F_2 : \wedge^2 H^{\bul}(\cL)  \to \cL
$$
of degree $-1$ such that
\begin{equation}
\label{F-1-F-2}
F_1([v_1, v_2]) - [F_1(v_1), F_1(v_2)] =
\pa F_2(v_1, v_2)\,,
\end{equation}
where $\pa$ is the differential on $\cL$\,.

It is convenient to enlarge the category of DG Lie algebras
by $\Linf$ algebras \cite{HS}, \cite{LS}. These more general algebras
have two important advantages. First, every sequence of quasi-isomorphisms
between $\Linf$ algebras $\cL$ and $\tcL$ can be shortened to
a single $\Linf$-quasi-isomorphism
$$
F: \cL \stackrel{\sim}{\to} \tcL\,.
$$
Second, for every DG Lie algebra (or possibly an $\Linf$ algebra)
$\cL$ there exists an $\Linf$-algebra structure on $H^{\bul}(\cL)$
such that $H^{\bul}(\cL)$ is quasi-isomorphic to $\cL$\,.

Thus if a DG Lie algebra $\cL$ is non-formal then this
$L_{\infty}$ algebra structure on $H^{\bul}(\cL)$ ``measures''
to what extent $\cL$ is far from being formal.

The most pedestrian way to introduce the notion of $\Linf$ algebra
is to start with the Chevalley-Eilenberg chain complex $C(\cL)$ of
a DG Lie algebra $\cL$\,.

As a graded vector space the Chevalley-Eilenberg chain complex is
the direct sum of all symmetric
powers of the desuspension $\bs^{-1} \cL$
of $\cL$
\begin{equation}
\label{CE}
C(\cL) = \sum_{k=1}^{\infty} \Big[ ( \bs^{-1}\cL )^{\otimes k}
\Big]^{S_k}\,.
\end{equation}
To introduce the differential we remark that $C(\cL)$ is
equipped with the following cocommutative comultiplication:
\begin{equation}
\label{copro-C}
\D\, :\, C(\cL)\mapsto
C(\cL) \otimes C(\cL)
\end{equation}
$$
\D(v)=0\,,
$$
\begin{equation}
\label{copro-eq}
\D (v_1, v_2, \dots , v_n)  =
\sum_{k=1}^{n-1} \sum_{\si \in {\rm Sh}(k,n-k)}
\pm (v_{\si(1)}, \dots, v_{\si(k)})
\otimes (
v_{\si(k+1)}, \dots , v_{\si(n)})\,,
\end{equation}
where $v_1, \dots,\, v_n$ are homogeneous elements
of $\bs^{-1}\cL$\,, ${\rm Sh}(k,n-k)$ is the set of $(k,n-k)$-shuffles
in $S_n$\,, and the signs are determined using the
Koszul rule.

To define the boundary operator $Q$ on (\ref{CE}) for
a DG Lie algebra $\cL$ we introduce the natural projection
\begin{equation}
\label{p}
p: C(\cL) \to \bs^{-1}\cL\,.
\end{equation}

It is not hard to see that if $Q$ is coderivation
of the coalgebra $C(\cL)$ in the sense of
the equation
$$
\D Q = (Q \otimes 1 + 1 \otimes Q) \D
$$
then $Q$ is uniquely determined by its composition
$p \circ Q$ with $p$\,.
This statement follows from the fact that $C(\cL)$
is a cofree\footnote{Strictly speaking $C(\cL)$ is a cofree cocommutative 
coalgebra without counit.} 
cocommutative coalgebra. The same statement
holds for cofree coalgebras of other types. (See Proposition 2.14
in \cite{GJ}.)

Thus we define the coboundary operator $Q$ in terms
of the differential $\pa$ and the Lie bracket $[\,,\,]$ by
requiring that it is a coderivation of the coalgebra
$C(\cL)$ and by setting
$$
p \circ Q(v) = - \pa v\,, \qquad
p \circ Q(v_1, v_2) = (-1)^{|v_1|+1} [v_1, v_2]\,,
$$
$$
p \circ Q(v_1, v_2, \dots, v_k) = 0\,, \qquad k > 2\,,
$$
where $v, v_1, \dots, v_k$ are homogeneous elements of $\cL$\,.

The equation $Q^2=0$ readily follows from the Leibniz rule
$$
\pa [v_1, v_2] = [\pa v_1, v_2] + (-1)^{|v_1|}[v_1, \pa v_2]\,,
$$
and the Jacobi identity:
$$
(-1)^{|v_1||v_3|}[[v_1, v_2], v_3] + {\rm c.p.}\{1,2,3\} = 0\,.
$$

To define a notion of $\Linf$ algebra we simply allow
the most general degree $1$ coderivation $Q$ of the coalgebra
$C(\cL)$ (\ref{CE}) satisfying the equation $Q^2=0$\,.
More precisely,
\begin{defi}
\label{defi-Linf}
An {\rm  $L_{\infty}$-algebra structure} on a graded vector space
$\cL$ is a degree $1$ coderivation $Q$ of the coalgebra $C(\cL)$ (\ref{CE})
with the comultiplication (\ref{copro-eq}) such that
$$
Q^2 = 0\,.
$$
Furthermore, an $\Linf$ {\rm morphism} $F$ from an $\Linf$ algebra
$(\cL, Q)$ to an $\Linf$ algebra $(\tcL, \tQ)$ is a
morphism of the corresponding DG cocommutative coalgebras:
\begin{equation}
\label{map-F}
F : (C(\cL), Q) \to  (C(\tcL), \tQ)\,.
\end{equation}
\end{defi}
The coderivation $Q$ is uniquely determined by the
degree $1$ maps
\begin{equation}
\label{Q-n}
Q_n =  p \circ Q \Big|_{(\bs^{-1}\cL)^{\otimes\, n}} :
(\bs^{-1}\cL)^{\otimes\, n} \to \tcL
\end{equation}
and we call them the structure maps of
the $\Linf$ algebra $(\cL, Q)$\,.

The equation $Q^2=0$ is equivalent to an
infinite collection of coherence relations
for $Q_n$'s. The first relation says that $Q^2_1=0$\,.
The second relation is the Leibniz identity for
$Q_1$ and $Q_2$\,. The third relation says that
the binary operation
$$
[\ga_1, \ga_2] = (-1)^{|\ga_1|+1} Q_2(\ga_1, \ga_2)
$$
satisfies the Jacobi identity up to homotopy and the
corresponding chain homotopy is exactly $Q_3$\,.
In particular, the $\Linf$ algebras with the zero higher
maps $Q_n$\,, $n>2$ are exactly the DG Lie algebras.

~\\
{\bf Remark.}
It is possible to define $\infty$ or homotopy versions
for algebras of other types. Although the intrinsic definition of
such $\infty$-versions
requires the language of operads \cite{BM}, \cite{Hinich-op}, \cite{GJ},
\cite{GK}, \cite{Markl}, we avoid this language here
for sake of accessibility and try to get by using the vague
analogy with the case of $L_{\infty}$ algebras.

~\\

By analogy with the coderivations every morphism
(\ref{map-F}) from the coalgebra  $(C(\cL), Q)$
to the coalgebra $(C(\tcL), \tQ)$
is uniquely determined by
its composition $p\circ F$ with the natural projection
$$
p: C(\tcL) \to \bs^{-1} \tcL\,.
$$
In other words, an $\Linf$-morphism is not a map
from $\cL$ to $\tcL$ but a collection of
(degree zero) maps:
\begin{equation}
\label{F-n}
F_n = p \circ F \Big|_{(\bs^{-1}\cL)^{\otimes n}} :
(\bs^{-1}\cL)^{\otimes n} \mapsto \bs^{-1} \tcL, \qquad n\ge 1
\end{equation}
compatible with the action of the symmetric groups $S_n$
and satisfying certain equations involving $Q$ and $\tQ$\,.

For this reason we reserve a special arrow $\leadsto$ for
$\Linf$-morphisms
\begin{equation}
\label{Linf}
F: \cL \leadsto \tcL\,.
\end{equation}
We call $F_n$ (\ref{F-n}) the
structure maps of the $\Linf$ morphism (\ref{Linf}).

It is not hard to see that the compatibility with
with $Q$ and $\tQ$ implies that $F_1$ in (\ref{F-n})
is a morphism between the cochain complexes $\bs^{-1}\cL$ and
$\bs^{-1}\tcL$\,.

\begin{defi}
An $\Linf$ {\rm quasi-isomorphism} from $\cL$ to $\tcL$
is an $\Linf$ morphism (\ref{Linf}) for which the
map $F_1: \bs^{-1}\cL \to \bs^{-1} \tcL$
is a quasi-isomorphism of cochain complexes.
\end{defi}

Following V. Hinich we have\footnote{V. Hinich \cite{Hinich} proved
this statement for $\infty$-versions of algebras over an
arbitrary quadratic Koszul operad.}
\begin{teo}[Lemma 4.2.1, \cite{Hinich}]
\label{Hinich-teo}
For every DG Lie algebra $\cL$ there exists
an $\Linf$ algebra structure $Q^H$ on $H^{\bul}(\cL)$
such that $\cL$ is quasi-isomorphic to
the $\Linf$-algebra
$$
(H^{\bul}(\cL), Q^H)\,.
$$
\end{teo}
The structure map $Q^H_1$ is zero and
$$
Q^H(\ga_1, \ga_2) = (-1)^{|\ga_1|+1}[\ga_1, \ga_2]\,,
$$
where $\ga_1, \ga_2\in H^{\bul}(\cL)$ and
$[\,,\,]$ is the induced Lie bracket on $H^{\bul}(\cL)$\,.

Thus, even if $\cL$ is a non-formal DG Lie algebra, the
Lie algebra structure on its cohomology $H^{\bul}(\cL)$ can be
corrected to an $\Linf$ algebra structure $Q^H$ such that
the $\Linf$ algebra $(H^{\bul}(\cL), Q^H)$ is quasi-isomorphic to
$\cL$\,.

The higher structure maps $Q^H_n$, $n >2$ depend on various choices
and formality of the DG Lie algebra (or more generally
$\Linf$ algebra) $(\cL, Q)$ means that these higher maps
can be chosen to be all zeros.

In the concluding section of this article we give an
example of a non-formal DG Lie algebra.

~\\
{\bf Remark.} Using the higher structure maps $Q^H_n$
one may construct operations which are independent of choices.
These operations are known as Massey-Lie products \cite{Babenko},
\cite{Halperin-Stasheff}, \cite{May-Massey}, \cite{Massey-Lie} and
formality of a DG Lie algebra $\cL$ implies that all the
Massey-Lie products are zero. A good
exposition on Massey-Lie products for the category of DG
commutative algebras is given in Section 2 of \cite{Babenko}.

\subsection{From $\Linf$ algebras back to DG Lie algebras}
As we see from Definition \ref{defi-Linf}
every $\Linf$ algebra $(\cL, Q)$ is defined by the DG cocommutative
coalgebra $C(\cL)$ with the codifferential $Q$\,. Using this
coalgebra we may construct a DG Lie algebra $\cR(\cL,Q)$ which is
quasi-isomorphic to the $\Linf$ algebra $(\cL,Q)$\,.
As a graded Lie algebra,
$\cR(\cL,Q)$ is the free Lie algebra generated by $C(\cL)$
\begin{equation}
\label{cR-cL}
\cR(\cL,Q) = \Lie\big(\, C(\cL)\, \big)\,.
\end{equation}
The differential on $\cR(\cL,Q)$ consists of two
parts. To define the first part we use, in the obvious way, the
codifferential $Q$\,. To define the second part we use the
comultiplication $\D$ (\ref{copro-eq}) on $C(\cL)$ viewing
(\ref{cR-cL}) as the dual version of the Harrison chain complex.

The construction of the free resolution $\cR(\cL,Q)$ for
an $\Linf$ algebra $(\cL, Q)$ is known in topology
as {\it the rectification} \cite{Board-V}. We describe this construction
in more details for a wider class of algebras in \cite{BLT}
(See Proposition 3 therein).

\subsection{DG Lie algebras and Maurer-Cartan elements}
\label{DG-MC}
Given a DG Lie algebra $\cL$ with the differential $\pa$ and
the Lie bracket $[\,,\,]$ over the field $\bbK$ we introduce
the DG Lie algebra $\cL[[\h]]$ over\footnote{Here $\h$ is a formal deformation
parameter.} the ring $\bbK[[\h]]$
extending $\pa$ and $[\,,\,]$ by $\bbK[[\h]]$-linearity.

\begin{defi}
A {\rm Maurer-Cartan (or MC) element} $\al$ of the DG Lie
algebra $\cL$ is a formal series $\al\in \h\, \cL^1[[\h]]$
of degree $1$ elements satisfying the equation
\begin{equation}
\label{MC-eq}
\pa \al + \frac{1}{2}[\al, \al] = 0\,.
\end{equation}
\end{defi}
MC elements may be formally compared to flat connections.

Let us consider the Lie algebra $\h\, \cL^0 [[\h]]$ of formal
series of degree zero elements in $\cL$\,. It is easy
to see that the Lie algebra
$\h\, \cL^0 [[\h]]$ is a projective limit of nilpotent Lie
algebras:
$$
\h\, \cL^0 [[\h]] \, \Big/ \, \h^N\, \cL^0 [[\h]]\,, \qquad N > 1\,.
$$
Therefore it can be exponentiated to the group
\begin{equation}
\label{mG}
\mG = \exp(\h\cL^0 [[\h]])\,.
\end{equation}

This group acts on the MC elements of $\cL$ according
to the formula:
\begin{equation}
\label{action}
\exp(\xi)\, \al =
\exp([\,\,,\xi])\, \al + f([\,\,,\xi])\, \pa \xi\,,
\end{equation}
where $f$ is the power series of the function
$$
f(x) = \frac{e^x - 1}{x}
$$
at the point $x=0$\,.
Two MC elements connected by the action
of the group $\mG$ may be thought of as equivalent
flat connections.

In this way we get the Goldman-Millson groupoid \cite{Goldman-M}
$\MC(\cL)$ which captures the formal one-parameter
deformation theory associated
to the DG Lie algebra $\cL$\,. Objects of this groupoid
are MC elements of $\cL$ and morphisms between two MC elements
$\al_1$ and $\al_2$ are elements of the group $\mG$ (\ref{mG})
which transform $\al_1$ to $\al_2$\,.

We denote by $\pi_0(\MC(\cL))$ the set of isomorphism
classes of the Goldman-Millson groupoid $\MC(\cL)$\,.

Every morphism  $\mu: \cL \to \tcL$ of DG Lie algebras gives us
an obvious functor
\begin{equation}
\label{mu-star}
\mu_* : \MC(\cL) \to \MC(\tcL)
\end{equation}
from the groupoid $\MC(\cL)$ to the groupoid $\MC(\tcL)$\,.

According to \cite{Ezra-Lie}, \cite{Goldman-M} and \cite{SS-MC}
we have the following theorem
\begin{teo}
\label{MC-teo}
If $\mu : \cL \to \tcL$ is a quasi-isomorphism of DG Lie
algebras then $\mu_*$ induces an bijection between
$\pi_0(\MC(\cL))$ and $\pi_0(\MC(\tcL))$\,.
\end{teo}

Theorem \ref{MC-teo} is an
immediate corollary of Proposition 4.9 in E. Getzler's
paper \cite{Ezra-Lie}.
According to W. Goldman, J. Millson,  M. Schlessinger and J. Stasheff,
\cite{Goldman-M}, \cite{SS-MC} every quasi-isomorphism
$\mu$ from $\cL$ to $\tcL$ induces an equivalence of groupoids
$\MC(\cL)$ and $\MC(\tcL)$ provided the DG Lie algebras $\cL$ and $\tcL$
are concentrated in non-negative degrees.

It turns out that if $\cL$ has elements in negative degrees then
$\MC(\cL)$ can be upgraded to a higher groupoid \cite{Ezra-higher}.
For $\Linf$ algebras this idea was
thoroughly developed by E. Getzler in \cite{Ezra-Lie}.
Then these results of E. Getzler were generalized by
A. Henriques in \cite{Henriques}.

\subsection{Twisting by a MC element}
\label{twisting}
Given a MC element $\al \in \h\,\cL[[\h]]$
of a DG Lie algebra $\cL$ we may modify the DG Lie algebra
structure on $\cL[[\h]]$ by switching to the new differential
\begin{equation}
\label{pa-twisted}
\pa + [\al, \,]\,.
\end{equation}
It is the MC equation (\ref{MC-eq}) which implies
the identity $(\pa + [\al, \,])^2 = 0$\,.

We denote the DG Lie algebra $\cL[[\h]]$ with the
differential (\ref{pa-twisted}) and the original Lie bracket
by $\cL^{\al}$
$$
\cL^{\al} = (\cL[[\h]], \pa + [\al, \,], [\,,\,])\,.
$$
Following D. Quillen \cite{Q} we call this procedure of
modifying the DG Lie algebra {\it twisting}.

It is obvious that every morphism $\mu: \cL \to \tcL$ of DG Lie
algebras extends by $\bbK[[\h]]$-linearity to the morphism
from $\cL^{\al}$ to $\tcL^{\mu(\al)}$\,. We denote this morphism
by $\mu^{\al}$.

We claim that
\begin{pred}[Proposition 1, \cite{thesis}]
\label{twist}
If $\mu:\cL\to \tcL$ is a quasi-isomorphism between DG Lie algebras
then so is the morphism
$$
\mu^{\al} : \cL^{\al} \to  \tcL^{\mu(\al)}\,.
$$
\end{pred}
For $\Linf$ algebras, the twisting procedure was
described in \cite{Ezra-Lie}.

\section{Algebraic structures on the Hochschild complexes}
\label{alg-str}
Let us introduce the algebras we are interested in.

First, we recall from \cite{G} that
\begin{defi}
\label{G-alg}
A graded vector space $V$ is a {\rm Gerstenhaber algebra}
if it is equipped with a graded commutative and associative
product $\wedge$ of degree $0$
 and a graded Lie bracket $[\,,\,]$ of degree $-1$\,.
These operations have to be compatible in the sense of the
following Leibniz rule
\begin{equation}
\label{cup-Lie}
[\ga, \ga_1 \wedge  \ga_2] =
[\ga, \ga_1] \wedge \ga_2 +
(-1)^{|\ga_1|(|\ga|+1)} \ga_1 \wedge [\ga, \ga_2]\,.
\end{equation}
\end{defi}

Second, we recall from \cite{Boris-book} that
\begin{defi}
\label{precalc}
A {\rm precalculus} is a pair of a Gerstenhaber algebra
$(V, \wedge, [,])$ and a graded vector
space $W$ together with
\begin{itemize}
\item a module structure
$i_{\bul} \,:\, V \otimes W
\mapsto W $ of the graded commutative algebra
$V$ on $W$\,,

\item an action $l_{\bul} \,:\, \bs^{-1} V\otimes W
\mapsto W$ of the Lie algebra
$\bs^{-1} V$ on
$W$ which is compatible with $i_{\bul}$ in
the sense of the following equations
\begin{equation}
\label{l-i}
i_a l_b - (-1)^{|a|(|b|+1)}
l_b i_a = i_{[a,b]}\,, \qquad
l_{a\wedge b} = l_a i_{b} + (-1)^{|a|}i_a l_b \,.
\end{equation}
\end{itemize}
\end{defi}
Furthermore,
\begin{defi}
\label{calc}
A {\rm calculus} is a precalculus
$(V,W, [,], \wedge, i_{\bul}, l_{\bul})$
with a degree $-1$ unary operation $\de$ on $W$
such that
\begin{equation}
\label{l-i-delta}
\delta \, i_{a} - (-1)^{|a|} i_{a} \, \delta =
l_a\,,
\end{equation}
and\footnote{Although $\de^2=0$\,, the operation $\de$ is
not considered as a part of the differential on $W$\,.}
$\de^2 =0$\,.
\end{defi}

The simplest examples of these algebraic structures
come from geometry. More precisely, if $M$ is smooth real
manifold then the graded vector space $T^{\bul}_{poly}(M)$
of polyvector fields on $M$ is a Gerstenhaber algebra.
The commutative product is simply the exterior product $\wedge$ and
the Lie bracket is the Schouten-Nijenhuis bracket
$[\,,\,]_{SN}$ \cite{Koszul}. This bracket is defined in the obvious
way for vector fields and for functions:
$$
[v_1, v_2]_{SN} = [v_1, v_2]\,, \qquad
[v, a]_{SN} = v(a)\,, \qquad [a,b]_{SN} = 0\,,
$$
$$
v, v_1, v_2\in T^{1}_{poly}(M)\,,
\qquad a,b \in T^{0}_{poly}(M) = C^{\infty}(M)\,.
$$
Then it is extended by Leibniz rule (\ref{cup-Lie})
to all polyvector fields.

Adding to polyvector fields the graded vector space
 $\Omb(M)$ of exterior forms with the reversed grading
\begin{equation}
\label{C-calculus}
(T^{\bul}_{poly}(M), \Omb(M))
\end{equation}
we get a calculus algebra.
The module structure on $\Omb(M)$ over the commutative algebra
$(T^{\bul}_{poly}(M), \wedge)$ is defined by the contraction
$$
i : T^{\bul}_{poly}(M) \otimes \Omb(M) \to \Omb(M)\,,
$$
the unary operation $\de$ is the de Rham differential $d$ and
the Lie algebra module structure on $\Omb(M)$ over
$(T^{\bul+1}_{poly}(M), [\,,\,]_{SN})$ is given by the Lie derivative
$$
l : T^{\bul+1}_{poly}(M) \otimes \Omb(M) \to \Omb(M)\,,
$$
\begin{equation}
\label{Lie-deriv}
l_{\ga} = d i_{\ga} - (-1)^{|\ga|} i_{\ga} d\,.
\end{equation}

Another example of calculus algebra comes from
noncommutative geometry. To give this example we start with
an arbitrary unital associative algebra $A$\,.

For every bimodule $U$ over the algebra $A$
we introduce the Hochschild cochain complex
\begin{equation}
\label{Cbu}
\Cbu(A, U) = \Hom (A^{\otimes \bul}, U)
\end{equation}
and the Hochschild chain complex
\begin{equation}
\label{Cbd}
\Cbd(A, U)= U \otimes A^{\otimes \bul}
\end{equation}
of $A$ with coefficients in $U$\,.

We reserve the same notation $\pa^{Hoch}$
both for the Hochschild
coboundary operator
\begin{equation}
\label{Hoch1}
(\pa^{Hoch} P)(a_0, a_1, \dots, a_k) = a_0 P(a_1, \dots, a_k)
-P(a_0 a_1, \dots, a_k) + P(a_0, a_1 a_2, a_3, \dots , a_k) -
\dots
\end{equation}
$$
+ (-1)^{k} P(a_0, \dots , a_{k-2}, a_{k-1} a_k)
+(-1)^{k+1} P(a_0, \dots , a_{k-2}, a_{k-1}) a_k
$$
on $\Cbu(A, U)$ and for Hochschild boundary operator
\begin{equation}
\label{Hoch11}
\pa^{Hoch}(u, a_1, \dots, a_m) =
(u a_1, a_2, \dots, a_m)
-(u, a_1 a_2, a_3, \dots,  a_m) + \dots
\end{equation}
$$
+(-1)^{m-1}(u, a_1, \dots, a_{m-2}, a_{m-1} a_m)
+(-1)^m (a_m u, a_1, a_2, \dots, a_{m-1})\,,
$$
$$
a_i\in A\,, \qquad  u\in U
$$
on $\Cbd(A, U)$\,.

We reserve the notation $HH^{\bul}(A,U)$ (resp. $HH_{\bul}(A,U)$)
for the Hochschild cohomology (resp. homology) groups of $A$
with coefficients in $U$\,.

In the case $U=A$ we simplify the notation for
the Hochschild complexes and for the (co)homology
groups:
\begin{equation}
\label{Cbu-A}
\Cbu(A) = \Cbu(A,A)\,,
\end{equation}
\begin{equation}
\label{Cbd-A}
\Cbd(A)= C_{-\bul}(A,A)\,,
\end{equation}
\begin{equation}
\label{HH-HH}
HH^{\bul}(A) = HH^{\bul}(A,A)\,,
\qquad
HH_{\bul}(A) = HH_{-\bul}(A,A)\,.
\end{equation}
For our purposes, we use the reversed
grading on the Hochschild chains of
$A$ with coefficients in $A$.

Here are the five algebraic operations
on the complexes $\Cbu(A)$ and $\Cbd(A)$ which play an
important role:
\begin{itemize}

\item the cup-product $\cup$
\begin{equation}
\label{cup}
P_1\cup P_2 (a_1, a_2, \dots, a_{k_1 + k_2}) =
P_1(a_1, \dots, a_{k_1}) P_2(a_{k_1+1}, \dots, a_{k_1+k_2})\,,
\end{equation}
$$
P_i \in C^{k_i}(A)\,,
$$

\item the Gerstenhaber bracket $[\,,\,]_G$
$$
[Q_1, Q_2]_{G} =
$$
\begin{equation}
\label{Gerst}
\sum_{i=0}^{k_1}(-1)^{i k_2}
Q_1(a_0,\,\dots , Q_2 (a_i,\,\dots,a_{i+k_2}),\, \dots,
a_{k_1+k_2}) -
(-1)^{k_1 k_2} (1 \leftrightarrow 2)\,,
\end{equation}
$$
Q_i \in  C^{k_i+1}(A)\,,
$$

\item the contraction $I_P$
of a Hochschild cochain $P \in C^k(A)$
with Hochschild chains
\begin{equation}
\label{I}
I_P (a_0, a_1, \dots, a_m) =
\begin{cases}
(a_0 P(a_1, \dots, a_k), a_{k+1}, \dots , a_m)\,, \qquad {\rm if} ~~ m\ge k \,, \\
0 \,, \qquad {\rm otherwise}\,,
\end{cases}
\end{equation}

\item the Lie derivative
of Hochschild chains along a Hochschild cochain $Q\in C^{k+1}(A)$
\begin{equation}
\label{L}
L_{Q}(a_0, a_1, \dots, a_m)=
\sum_{i=0}^{m-k}(-1)^{ki} (a_0 , \dots,
Q(a_i, \dots, a_{i+k}), \dots, a_m) +
\end{equation}
$$
\sum_{j=m-k}^{m-1}(-1)^{m(j+1)} ( Q(a_{j+1}, \dots, a_m, a_0, \dots, a_{k+j-m}),
a_{k+j+1-m}, \dots , a_j )\,,
$$

\item and Connes' operator $B: \Cbd(A) \to C_{\bul-1}(A)$
\begin{equation}
\label{B}
\begin{array}{c}
\displaystyle
B (a_0, \dots, a_m) =
\sum_{i=0}^m \Big( (-1)^{mi}
(1, a_i, \dots, a_m, a_0, \dots, a_{i-1}) \\[0.3cm]
\displaystyle
+ (-1)^{m i} (a_{i}, 1, a_{i+1}, \dots, a_m,
a_0, \dots, a_{i-1}) \Big)\,.
\end{array}
\end{equation}

\end{itemize}
All these operations are compatible with the differential
$\pa^{Hoch}$ (\ref{Hoch1}), (\ref{Hoch11}).
Therefore they induce the corresponding operations
on the level of cohomology.

There are several identities involving the operations
$[\,,\,]_G$, $L$, and $B$.

First, the Gerstenhaber
bracket $[\,,\,]_G$ is a Lie bracket on $\sCbu(A)$ and hence
$\sCbu(A)$ is a DG Lie algebra. The operation $L$ (\ref{L})
gives us an action of the DG Lie algebra $\sCbu(A)$ on
Hochschild chains. In other words,
\begin{equation}
\label{L-brack}
L_{Q_1} L_{Q_2} - (-1)^{(|Q_1|+1)(|Q_2|+1)} L_{Q_2} L_{Q_1}
= L_{[Q_1, Q_2]_G}
\end{equation}
and hence $\Cbd(A)$ is a DG Lie algebra module over
the DG Lie algebra $\sCbu(A)$\,.

The Connes cyclic operator $B$ (\ref{B}) is used in the
definitions of different variants of cyclic chain complex
\cite{Boris-book}, \cite{Loday}. All these variants have
the form\footnote{This notation is due to E. Getzler.}
\begin{equation}
\label{CC}
CC^{\cW}_{\bul}(A)  = (\Cbd(A)[[u]]\, \otimes_{\bbK[u]}\, \cW, \pa^{Hoch} + uB)\,,
\end{equation}
where $u$ is an auxiliary variable of degree $2$ and
$\cW$ is a $\bbK[u]$-module.

We are interested in two particular cases:

--- if $\cW= \bbK[u]$ then $CC^{\cW}_{\bul}(A)$ is
called the negative cyclic complex
\begin{equation}
\label{CC-neg}
CC^-_{\bul}(A)  = (\Cbd(A)[[u]], \pa^{Hoch} + uB)\,,
\end{equation}

--- if $\cW= \bbK$ with $u$ acting by zero then
$CC^{\cW}_{\bul}(A)$ is nothing but the Hochschild chain
complex (\ref{Cbd-A}) of $A$\,.

Since the Connes cyclic operator $B$
(\ref{B}) is compatible with the ``Lie derivative'' $L$
(\ref{L}) in the sense of the equation
$$
B L_P - (-1)^{|P|+1} L_P B = 0\,, \qquad
P \in \Cbu(A)
$$
any variant of the cyclic chain complex (\ref{CC}) is a DG module over
the DG Lie algebra $\sCbu(A)$\,.

Unfortunately, the cup-product $\cup$ and the Gerstenhaber
bracket $[\,,\,]_G$ do not satisfy the Leibniz rule.
So $\cup$ and $[\,,\,]_G$ do not give us a Gerstenhaber
algebra structure on $\Cbu(A)$\,. Similarly, the operations
$\cup$, $[\,,\,]_G$, $I$, $L$, $B$ do not give us
a calculus algebra on the pair $(\Cbu(A), \Cbd(A))$\,.
However, the
required identities hold on the level of cohomology and
we have
\begin{pred}[M. Gerstenhaber \cite{G}]
\label{HH-Gerst}
The cup-product (\ref{cup}) and the bracket
(\ref{Gerst}) induce on $HH^{\bul}(A)$ a structure
of a Gerstenhaber algebra.
\end{pred}
and
\begin{pred}[Yu. Daletski, I. Gelfand, and B. Tsygan
\cite{DGT}]
\label{HH-calc}
The operations (\ref{cup}), (\ref{Gerst}),
(\ref{I}), (\ref{L}), and (\ref{B}) induce on the pair
$$
(HH^{\bul}(A), HH_{\bul}(A))
$$
a structure of a calculus.
\end{pred}

Let us also recall the Van den Bergh
duality theorem:
\begin{teo}[M. Van den Bergh, \cite{VB}]
\label{VdB}
If $A$ is a finitely generated bimodule
coherent\footnote{An algebra $A$ is called bimodule coherent
if every map between finite rank free $A$-bimodules
has a finitely generated kernel (see Definition 3.5.1
in \cite{Vitya}).}
algebra of finite Hochschild dimension $d$\,,
\begin{equation}
\label{Ext}
HH^m (A, A \otimes A) =
\begin{cases}
V_A  & {\rm if} \quad m = d\,,\\
0  & {\rm otherwise}\,,
\end{cases}
\end{equation}
where $V_A$ is an invertible\footnote{A $A$-bimodule $V$ is
called invertible if there is a $A$-bimodule $\tV$ such
that $V \otimes_A \, \tV \cong A$\,.}
$A$-bimodule then
for every $A$-bimodule $U$
$$
HH^{\bul}(A, U) \cong HH_{d - \bul}(A, V_A \otimes_A U)\,.
$$
\end{teo}
In Equation (\ref{Ext}) $A\otimes A$ is considered
as a bimodule over $A$ with respect
to the external $A$-bimodule structure.
It is the internal $A$-bimodule
structure which equips all the cohomology groups
$HH^{\bul} (A, A \otimes A)$
with a structure of $A$-bimodule.

We refer to $V_A$ as {\it the Van den Bergh dualizing
module} of $A$\,.

\section{Formality theorems}
\label{F}
The famous Kontsevich's formality theorem can be
formulated as
\begin{teo}[M. Kontsevich, \cite{K}]
\label{Maxim}
Let $M$ be a smooth real manifold and $A = C^{\infty}(M)$
be the algebra of smooth functions on $M$\,. Then the DG Lie algebra
$(\sCbu(A), \pa^{Hoch},[\,,\,]_G)$ of Hochschild cochains
is quasi-isomorphic to the graded Lie algebra
$(T^{\bul+1}_{poly}(M), [\,,\,]_{SN})$ of polyvector fields.
\end{teo}
{\bf Remark.} Since $A= C^{\infty}(M)$ is a topological algebra the definition
of Hochschild cochains for $A$ requires some precaution.
By the Hochschild cochains of the algebra $A= C^{\infty}(M)$ we
mean the polydifferential operators on $M$\,.

To prove Theorem \ref{Maxim} M. Kontsevich gave \cite{K} an
explicit construction of an $\Linf$ quasi-isomorphism
\begin{equation}
\label{cK}
\cK : T^{\bul+1}_{poly}(\bbR^d) \leadsto \sCbu(C^{\infty}(\bbR^d))
\end{equation}
from the graded Lie algebra $T^{\bul+1}_{poly}(\bbR^d)$ of polyvector fields
to the DG Lie algebra $\sCbu(C^{\infty}(\bbR^d))$ of polydifferential
operators on $\bbR^d$\,. This construction involves very
interesting integrals over compactified configuration spaces of
points on the upper half plane. Tedious questions about
choices of signs were thoroughly addressed in paper \cite{signs}.

In order to extend this result to an arbitrary smooth
manifold M. Kontsevich used what is called
the Gelfand-Fuchs
trick \cite{GF} or the formal geometry \cite{G-Kazh}
in the sense of I.M. Gelfand and D.A. Kazhdan.

This step of globalization was discussed later
in more details by A. Cattaneo, G. Felder and Tomassini
in \cite{CFT}, by M. Kontsevich in Appendix 3 in \cite{K-alg},
by the first author in \cite{CEFT},
by A. Yekutieli in \cite{Ye} and by M. Van den Bergh
in \cite{VdB-glob}.

To describe the first immediate corollary of Theorem
\ref{Maxim} we recall from \cite{Bayen} and \cite{Ber} that a star-product
on the manifold
$M$ is an $\bbR[[\h]]$-linear associative product
on $C^{\infty}(M)[[\h]]$ of the form
\begin{equation}
\label{star}
a * b = a b + \sum_{k=1}^{\infty} \h^k \Pi_k(a,b)\,,
\end{equation}
$$
a,b \in C^{\infty}(M)[[\h]]\,,
$$
where $\Pi_k$ are bidifferential operators.

Since $a*b = a b ~~ {\rm mod}~~ \h$ the star-product
(\ref{star}) should be viewed as an associative
and not necessarily commutative formal deformation
of the ordinary product of functions on $M$\,.

Two star-products $*$ and $\tilde{*}$ are called
equivalent if there exists a formal series of differential
operators
\begin{equation}
\label{T}
T = Id + \h T_1 + \h^2 T_2 + \dots
\end{equation}
which starts from the identity and intertwines the star-products
$*$ and $\tilde{*}$:
\begin{equation}
\label{equiv}
T(a * b) = T(a) \,\tilde{*}\, T(b)\,.
\end{equation}

It easy to see that the associativity property of the star-product
(\ref{star}) is equivalent to the MC equation
$$
\pa^{Hoch} \Pi + \frac{1}{2}[\Pi,\Pi]_G = 0
$$
for the element
\begin{equation}
\label{Pi}
\Pi = \sum_{k=1}^{\infty} \h^k \Pi_k\,.
\end{equation}
Thus MC elements of the DG Lie algebra $\sCbu(\,C^{\infty}(M)\,)$
are exactly the star-products on $M$\,.

Furthermore, it is not hard to identify the intertwiners between
star-products with morphisms of the Goldman-Millson groupoid
$\MC(\sCbu(\,C^{\infty}(M)\,))$\,.

For the graded Lie algebra $T^{\bul+1}_{poly}(M)$ of polyvector fields
MC elements are the formal Poisson structures. These are
formal power series of bivectors
\begin{equation}
\label{pi}
\pi = \sum_{k=1}^{\infty} \h^k \pi_k\,, \qquad
\pi_k \in T^{2}_{poly}(M)
\end{equation}
satisfying the Jacobi relation
\begin{equation}
\label{Jac}
[\pi, \pi]_{SN} = 0\,.
\end{equation}
We say that two formal Poisson structures $\pi$ and
$\tpi$ are equivalent if they are connected by
the adjoint action of the group
\begin{equation}
\label{mG-for-Tpoly}
\exp\Big(\h\, T^1_{poly}(M)[[\h]] \Big)\,.
\end{equation}

Thus combining Theorem \ref{MC-teo} and Theorem \ref{Maxim}
we get the following corollary
\begin{cor}
\label{classif-def}
Equivalence classes of star products are in
bijection with the equivalence classes of
formal Poisson structures.
\end{cor}

Given a star-product $*$ on $M$, we call the
corresponding equivalence class of formal Poisson structures
{\it Kontsevich's class} of the star-product.
Notice that the action of the group (\ref{mG-for-Tpoly})
does not change the first term $\pi_1$ of the series
(\ref{pi}). Thus $\pi_1$ does not depend on the
choice of the representative (\ref{pi}) of Kontsevich's
class. We refer to $\pi_1$ as the Poisson bivector corresponding
to the star-product $*$.

Twisting the DG Lie algebra $\sCbu(\,C^{\infty}(M)\,)$ by
the MC element $\Pi$ (\ref{Pi}) corresponding to the star-product
$*$ (\ref{star})
we get the DG Lie algebra $\sCbu(\,A_{\h}\,)$ of Hochschild cochains
for the deformation quantization algebra
$A_{\h}= (C^{\infty}(M)[[\h]], *)$\,.

Twisting the graded Lie algebra of polyvector fields
$T^{\bul+1}_{poly}(M)$ by the MC element $\pi$ (\ref{pi})
we get the Poisson cochain complex \cite{Lich}
\begin{equation}
\label{HP-up}
(T^{\bul+1}_{poly}(M)[[\h]], [\pi,\,]_{SN})
\end{equation}
of $\pi$\,.
Omitting the shift we refer to the cohomology of the complex (\ref{HP-up})
as the Poisson cohomology $HP^{\bul}(M, \pi)$ of the formal Poisson
structure $\pi$
\begin{equation}
\label{HP-up1}
HP^{\bul}(M,\pi) := H^{\bul}(T^{\bul}_{poly}(M)[[\h]], [\pi,\,]_{SN})\,.
\end{equation}

Thus, combining Proposition \ref{twist} with Theorems \ref{MC-teo} and
\ref{Maxim}, we get the following corollary
\begin{cor}
\label{HH-HP-up}
If $*$ is a star-product whose Kontsevich's class
is represented by the formal Poisson structure $\pi$
then the Hochschild cohomology of the deformation quantization
algebra $A_{\h}= (C^{\infty}(M)[[\h]], *)$ is isomorphic
to Poisson cohomology $HP^{\bul}(M, \pi)$ of the
formal Poisson structure $\pi$\,.
\end{cor}

\subsection{Alternative approach to Theorem \ref{Maxim}}
As we mentioned above, the operations $\cup$ (\ref{cup}) 
and $[\,,\,]_G$
(\ref{Gerst}) do not equip the Hochschild cochain
complex $\Cbu(A)$ with a Gerstenhaber algebra structure.

However, using the solution of Deligne's Hochschild cohomological
conjecture  \cite{BF}, \cite{K-Soi}, \cite{M-Smith},
\cite{Dima-DG}, \cite{Sasha-Proof} one can show that
the operations $\cup$ and $[\,,\,]_G$ can be upgraded to
an $\infty$-Gerstenhaber algebra whose multiplications are
expressed in terms of the cup-product $\cup$ and
insertions of cochains into a cochain.

This $\infty$-Gerstenhaber structure depends on
the choice of Drinfeld's associator \cite{Drinfeld}.
At this moment it is not known whether the homotopy
type of this $\infty$-Gerstenhaber structure depends on
this choice.
However,
\begin{teo}[Theorem 2.1, \cite{Dima-Proof}]
\label{non-spoiled}
The $\Linf$ algebra part of the
$\infty$-Gerstenhaber structure on
Hochschild cochains coincides with
the DG Lie algebra structure
given by the Hochschild differential
and the Gerstenhaber bracket (\ref{Gerst}).
\end{teo}

In 1998 the second author proposed a completely different
proof \cite{Dima-Proof} of Theorem \ref{Maxim} in the case
when $A$ is the polynomial algebra over an arbitrary field $\bbK$
of characteristic zero. We would like to refer the reader to
excellent V. Hinich's  exposition \cite{Hinich} of this approach.

In \cite{Dima-Proof} it was shown that the Gerstenhaber algebra
$T^{\bul}_{poly}(\bbK^d)$
of polyvector fields on the affine space $\bbK^d$ is {\it intrinsically formal}.
In other words, there is no room for cohomological obstructions
to the formality of any $\infty$-Gerstenhaber algebra whose cohomology
is the Gerstenhaber algebra $T^{\bul}_{poly}(\bbK^d)$\,.

Therefore, the Hochschild cochain complex $\Cbu(\bbK[x_1, \dots, x_d])$
with the $\infty$-Gerstenhaber structure coming from the
solution of Deligne's conjecture is formal.

Combining this observation with Theorem \ref{non-spoiled} we
immediately deduce Theorem \ref{Maxim}.

At this moment there is a more explicit proof of the
formality for the $\infty$-Gerstenhaber structure for
a wider class of algebras. More precisely,
\begin{teo}[\cite{BLT}]
\label{BLT-teo}
For every regular commutative algebra $A$ over
a field $\bbK$ of characteristic zero the $\infty$-Gerstenhaber
algebra $\Cbu(A)$ of Hochschild cochains is formal.
\end{teo}
An analogous statement in the Lie algebroid setting
was proved by D. Calaque and M. Van den Bergh in \cite{CV-Ginfty}.

\subsection{Formality theorems for Hochschild and cyclic chains}
Hochschild chains (\ref{Cbd-A}) enter this picture in
a very natural way
\begin{teo}
\label{chains}
Let $M$ be a smooth real manifold and $A = C^{\infty}(M)$
be the algebra of smooth functions on $M$\,. Then the
pair ``DG Lie algebra $\sCbu(A)$ and its DG module $\Cbd(A)$''
is quasi-isomorphic to the pair ``graded Lie algebra
$T^{\bul+1}_{poly}(M)$ and its module $\Omb(M)$''.
\end{teo}
This statement was formulated as a conjecture by
the third author in \cite{Tsygan}.
It was proved in \cite{Sh} by B. Shoikhet for the case $M=\bbR^d$\,.
The step of globalization was performed in the thesis of the
first author \cite{FTHC}, \cite{thesis}.

As well as Theorem \ref{Maxim}, Theorem \ref{chains}
requires a technical amendment. More precisely, the
space $C_k(A)$ of degree $k$ Hochschild chains for
the algebra $A = C^{\infty}(M)$ should be replaced by
the space
of $\infty$-jets near the main diagonal of
the product $M^{\times\, k}$\,.

In \cite{Tsygan} the third author also conjectured the
formality of the cyclic complexes (\ref{CC}) as
DG Lie algebra modules over $\sCbu(A)$\,.
In paper \cite{TT} the second author and the third author
proposed a plan on how this cyclic
conjecture can be proved.

However, in \cite{W} T. Willwacher showed elegantly that
Shoikhet's $\Linf$ quasi-isomorphism \cite{Sh} is compatible
with the Connes cyclic operator (\ref{B}). This observation readily
settled in the positive Tsygan's cyclic formality conjecture:
\begin{teo}[T. Willwacher, \cite{W}]
\label{cyclic}
Let $M$ be a smooth real manifold and $A = C^{\infty}(M)$
be the algebra of smooth functions on $M$\,. If $\cW$ is a
$\bbR[u]$-module of finite projective dimension then
the pair
``DG Lie algebra $\sCbu(A)$ and its DG module $CC^{\cW}_{\bul}(A)$''
is quasi-isomorphic to the pair ``graded Lie algebra
$T^{\bul+1}_{poly}(M)$ and its DG module
$(\Omb(M)[[u]]\, \otimes_{\bbR[u]}\, \cW\,,\, u\, d)$''.
\end{teo}
If we set $\cW= \bbR$ with $u$ acting by zero then
$CC^{\cW}_{\bul}(A)$ turns to the Hochschild chain complex
and the DG module
$$
(\Omb(M)[[u]]\, \otimes_{\bbR[u]}\, \cW\,,\, u\, d)
$$
turns to the module $\Omb(M)$ with the zero differential.
Thus Theorem \ref{chains} is a corollary of Theorem \ref{cyclic}.

Given a Lie algebra $V$ and its module $W$\,, we can form the semi-direct
product $V\oplus W$ in which $W$ is an Abelian Lie algebra.
It is clear that pairs ``Lie algebra $V$ and its module $W$''
can be identified with such semi-direct products.

Using this idea one can generalize the twisting procedure
we described in Subsection \ref{twisting} to DG Lie algebra
modules.

Thus, if $*$ is a star-product on $M$ then twisting
the DG Lie algebra module of negative cyclic chains (\ref{CC-neg})
of the algebra  $A = C^{\infty}(M)$ by the corresponding
MC element $\Pi$ (\ref{Pi}) we get the negative
cyclic complex
$$
CC^{-}_{\bul}(A_{\h})
$$
for the deformation quantization algebra
$A_{\h} = (C^{\infty}(M)[[\h]], *)$\,.

Let $\pi$ be a formal Poisson structure (\ref{pi})\,.
Regarding $\pi$ as a MC element of the graded Lie algebra
$T^{\bul}_{poly}(M)$ and twisting its DG Lie algebra module
$$
(\Omb(M)[[u]]\, , \, u\, d)
$$
we get the DG Lie algebra module
\begin{equation}
\label{formy-l-ud}
(\Omb(M)[[u]][[\h]]\,,\, l_{\pi} + u\, d)
\end{equation}
over the DG Lie algebra $(T^{\bul}_{poly}(M)[[\h]], [\pi\,,\,]_{SN} )$
with the Lichnerowicz differential $[\pi\,,\,]_{SN}$\,.

Generalizing Proposition \ref{twist} to DG Lie algebra
modules in the obvious way we get the following corollary
of Theorem \ref{cyclic}
\begin{cor}
\label{cyclic-comp}
If $*$ is a star-product on $M$ whose Kontsevich's class
is represented by the formal Poisson structure $\pi$ then
the complex (\ref{formy-l-ud}) computes the negative cyclic
homology
$$
HC^{-}_{\bul}(A_{\h}) = H^{\bul}(\, CC^{-}_{\bul}(A_{\h})\, )
$$
of the deformation quantization algebra $A_{\h} = (C^{\infty}(M)[[\h]], *)$\,.
\end{cor}
{\bf Remark.} Applying similar arguments to the Hochschild chain
complex and to the exterior forms we get an isomorphism
between Hochschild homology of $\hA$ and the Poisson homology
\cite{JLB}, \cite{Koszul} of $\pi$\,.

\subsection{Formality of the $\infty$-calculus algebra
$(\Cbu(A), \Cbd(A))$}
To include the operations $\cup$ (\ref{cup}),  $[\,,\,]_G$ (\ref{Gerst}),
$I$ (\ref{I}), $L$ (\ref{L}), and $B$ (\ref{B}) on
the pair
\begin{equation}
\label{para}
(\Cbu(A), \Cbd(A))
\end{equation}
into the picture we need to find a correct algebraic
structure on (\ref{para}).

As we already mentioned above the operations
$\cup$ (\ref{cup}),  $[\,,\,]_G$ (\ref{Gerst}),
$I$ (\ref{I}), $L$ (\ref{L}), and $B$ (\ref{B}) satisfy
the identities of the calculus algebra only up to homotopy.
So a calculus algebra is not a correct algebraic structure
for this situation.

Luckily the Kontsevich-Soibelman solution \cite{K-Soi1} of
the chain version of Deligne's conjecture implies that
the operations $\cup$,  $[\,,\,]_G$, $I$, $L$, and $B$
can be upgraded to an $\infty$-calculus algebra whose multiplications
are expressed in terms of the cup-product $\cup$\,,
insertions of cochains into a cochain, and insertions of components of
a chain into cochains which respect the cyclic order on these
components.

In \cite{Cepochki} we show that
\begin{teo}[Corollary 4, \cite{Cepochki}]
\label{calc-formality}
For every smooth real manifold $M$ the
$\infty$-calculus algebra
\begin{equation}
\label{para-M}
\big( \Cbu(C^{\infty}(M)), \Cbd(C^{\infty}(M)) \big)
\end{equation}
is quasi-isomorphic to the calculus algebra
\begin{equation}
\label{para-T-Om}
(T^{poly}_{\bul}(M), \Omb(M))
\end{equation}
of polyvector fields and exterior forms on $M$\,.
\end{teo}
We prove this theorem using our construction from
$\cite{BLT}$ for the $\infty$-Gerstenhaber algebra
on $\Cbu(C^{\infty}(M)$
and the Morita equivalence between the algebra
of differential operators on exterior forms and the
algebra of differential operators on functions.
Unlike in \cite{BLT}, we did not produce an explicit
sequence of quasi-isomorphisms connecting (\ref{para-M})
and (\ref{para-T-Om}). So Theorem \ref{calc-formality}
has a status of an existence theorem.

Theorem 4 from \cite{Cepochki} implies that the algebra
structure on (\ref{para}) given by operations $[\,,\,]_G$, $L$ and
$B$ do not have higher homotopy corrections inside the
$\infty$-calculus structure. Thus the cyclic formality theorem
of T. Willwacher (Theorem \ref{cyclic}) is a corollary
of Theorem \ref{calc-formality}.

We should remark that Kontsevich's $\Linf$ quasi-isomorphism $\cK$ (\ref{cK})
has a subtle compatibility property with the cup-product.
To formulate this property we set $\pi$ to be a formal Poisson structure (\ref{pi})
on $\bbR^d$ and $\hA = (C^{\infty}(\bbR^d)[[\h]], *)$ be a deformation
quantization algebra whose Kontsevich's
class is represented by $\pi$\,.
Twisting  Kontsevich's $\Linf$ quasi-isomorphism $\cK$ (\ref{cK})
by $\pi$ we get the $\Linf$ quasi-isomorphism
\begin{equation}
\label{cK-tw}
\cK^{\pi} : \big(T^{\bul+1}_{poly}(\bbR^d)[[\h]],  [\pi,\,]_{SN} \big)
\leadsto \sCbu(\hA)
\end{equation}
which, in turn, induces an isomorphism from
the Poisson cohomology
\begin{equation}
\label{HP-up-Rd}
HP^{\bul}(\bbR^d,\pi) := H^{\bul}(T^{\bul}_{poly}(\bbR^d)[[\h]], [\pi,\,]_{SN})\,.
\end{equation}
to the Hochschild cohomology
\begin{equation}
\label{HH-up-Rd}
HH^{\bul}(\hA)
\end{equation}
of $\hA$\,.

The exterior product $\wedge$ turns the Poisson cohomology (\ref{HP-up-Rd})
into a graded commutative algebra. Similarly, the cup-product
$\cup$ (\ref{cup}) turns the Hochschild cohomology (\ref{HH-up-Rd})
into a graded commutative algebra. Due to \cite{K} and
\cite{Torossian-cup} we have the following theorem
\begin{teo}
\label{with-cup}
The isomorphism from the Poisson cohomology (\ref{HP-up-Rd})
to the Hochschild cohomology (\ref{HH-up-Rd}) induced by
the $\Linf$ quasi-isomorphism (\ref{cK-tw}) is an
isomorphism of graded commutative algebras.
\end{teo}
This property was used in \cite{CV} by D. Calaque and M. Van den
Bergh to prove C\u{a}ld\u{a}raru's conjecture \cite{Cald} on
Hochschild structure of an algebraic variety.

Shoikhet's quasi-isomorphism of $\Linf$-modules \cite{Sh} has
a similar subtle compatibility property with the contraction $I$ (\ref{I}).
This property is proved by D. Calaque and C. Rossi in \cite{CR-cap}.
In paper \cite{CR-cap} the authors also applied this result
to obtain a version of the Duflo isomorphism on coinvariants.

These results indicate that there should be a bridge between
Kontsevich's construction \cite{K} (resp. Shoikhet's construction
\cite{Sh}) and the construction
of the second author \cite{Dima-Proof} (resp. the construction
in \cite{Cepochki}).

\subsection{Formality theorems for Hochschild and cyclic complexes
in the algebraic geometry setting}

In the algebraic geometry setting Theorem \ref{Maxim} has
a formulation which does not need the amendment about
the nature of cochains. In other words, we do not need
to restrict ourselves to the polydifferential operators.
\begin{teo}
\label{Maxim-alg}
Let $X$ be a smooth affine variety over a field $\bbK$ of
characteristic zero and $A=\cO_X(X)$ be the algebra of regular functions
on $X$\,. Then the DG Lie algebra
$(\sCbu(A), \pa^{Hoch},[\,,\,]_G)$ of Hochschild cochains
is quasi-isomorphic to the graded Lie algebra
$(\wedge^{\bul+1}_A \, \Der(A), [\,,\,]_{SN})$ of polyderivations
of $A$\,.
\end{teo}
This statement is an immediate corollary of
Theorem 2 and Theorem 4 from \cite{BLT}. It can be also
extracted from M. Kontsevich's paper \cite{K-alg} on deformation
quantization of algebraic varieties. However, Kontsevich's approach
requires that the base field $\bbK$ contains reals.

Beyond the affine case
it no longer makes sense to talk about global sections.
Thus we need to reformulate the question for the sheaves
of Hochschild cochains.

According to R. Swan \cite{Swan} and A. Yekutieli \cite{Y}
an appropriate candidate for the sheaf of Hochschild cochains
on an arbitrary smooth algebraic variety $X$ is the
sheaf of polydifferential operators with regular coefficients.
We denote this sheaf by $\Cbu(\cO_X)$\,. The Gerstenhaber
bracket (\ref{Gerst}) equips $\sCbu(\cO_X)$ with a structure of
a sheaf of DG Lie algebras.

Due to \cite{BLT}, \cite{VdB-glob}, and \cite{Ye} we have
\begin{teo}
\label{Maxim-nonaffine}
For every smooth algebraic variety $X$ over
a field $\bbK$ of characteristic zero the sheaf
of DG Lie algebra $\sCbu(\cO_X)$ is quasi-isomorphic
to the sheaf $\wedge^{\bul+1}T_X$ of polyvector fields
with the Schouten-Nijenhuis bracket.
\end{teo}
For applications of this theorem to deformation
quantization in the setting of algebraic geometry
we refer the reader to papers \cite{CH}, \cite{K-alg},
and \cite{Ye}.

Hochschild chains can also be added into this picture.
An appropriate candidate for the sheaf of Hochschild chains
on an algebraic variety $X$ is the sheaf of polyjets.
\begin{equation}
\label{polyjets}
\Cbd(\cO_X) = \cHom_{\cO_X}(C^{-\bul}(\cO_X), \cO_X)\,,
\end{equation}
where $\cHom$ denotes the sheaf-Hom and $\Cbu(\cO_X)$
is considered with its natural left $\cO_X$-module
structure.

Using the isomorphism between $X$ and the main diagonal
of the product $X^{\times \, k}$ we may identify local
sections of $C_k(\cO_X)$ with $\infty$-jets on $X^{\times \, k}$
near its main diagonal.

To introduce sheaves of cyclic chains we introduce
an auxiliary variable $u$ of degree $2$ and consider
a $\bbK[u]$-module $\cW$ as a constant sheaf on $X$\,.
Then to every such module we assign a sheaf
of cyclic chains
\begin{equation}
\label{CC-sheaf}
CC^{\cW}_{\bul}(\cO_X)  = (\Cbd(\cO_X)[[u]]\, \otimes_{\bbK[u]}\, \cW,
\pa^{Hoch} + uB)\,,
\end{equation}

The operation $L$ (\ref{L}) equips $CC^{\cW}_{\bul}(\cO_X)$
with a structure of a sheaf of DG Lie algebra
modules over the sheaf of DG Lie algebras $\Cbu(\cO_X)$\,.

Theorem 4 and 5 from \cite{Cepochki} implies
the following statement
\begin{teo}
\label{cyclic-alg}
Let $X$ be a smooth algebraic variety over
a field $\bbK$ of characteristic zero.
If $\cW$ is a
$\bbR[u]$-module of finite projective dimension then
the pair
``the sheaf of DG Lie algebras $\sCbu(\cO_X)$ and
the sheaf of its DG modules $CC^{\cW}_{\bul}(\cO_X)$''
is quasi-isomorphic to the pair ``the sheaf of graded Lie algebras
$\wedge^{\bul+1} T_X$ and the sheaf its DG modules
$(\Omb_X[[u]]\, \otimes_{\bbR[u]}\, \cW\,,\, u\, d)$''.
\end{teo}
{\bf Remark.} Using the construction of M. Van den Bergh
\cite{VdB-glob} and A. Yekutieli \cite{Y} as well as the
results of B. Shoikhet \cite{Sh} and T. Willwacher \cite{W}
it is possible to prove Theorem \ref{cyclic-alg}
under the assumption that the base field $\bbK$ contains reals.

\section{More application of formality theorems}
\label{more-applic}
Due to Corollary \ref{classif-def} the deformation
class of a star product $*$ is uniquely determined by
Kontsevich's class which is the equivalence class of a
formal Poisson structure. For this reason Kontsevich's class
is, sometimes, referred to as the characteristic class of
a star-product.

In many cases deformation quantization algebras are
not obtained via formality theorems.
In such situations it may be hard to find
Kontsevich's class of the
deformation quantization algebra.

However, it is often possible to extract some information
about representatives of this class using
the homological properties of the deformation
quantization algebra.

An example of this situation is provided by the unimodularity criterion
from \cite{Modular}.

To describe this criterion we set $X$ to be a smooth affine
variety with the trivial canonical bundle.
The triviality of the canonical bundle implies that there
exists a nowhere vanishing top degree exterior form $\Vol$
on $X$\,.

A formal Poisson structure $\pi$ (\ref{pi}) is called
{\it unimodular} \cite{BZ}, \cite{Alan} if there exists a formal
power series
$$
\Vol_{\h} = \Vol + \sum_{k=1}^{\infty} \h^k \Vol_k
$$
of top degree exterior forms starting with a nowhere
vanishing form $\Vol$ and such that
\begin{equation}
\label{unimod}
l_{\pi} \Vol_{\h}  = 0\,,
\end{equation}
where $l_{\pi}$ is the Lie derivative (\ref{Lie-deriv}).

Let $\hA$ be a deformation quantization algebra
$(A[[\h]], *)$ of $X$\,. It turns out that $\hA$ has the
Hochschild dimension equal to the dimension of $X$\,.
Furthermore, $\hA$ satisfies all the conditions of
Theorem \ref{VdB}.  Due to \cite{Modular}
we have the following
homological unimodularity criterion:
\begin{teo}[Theorem 3, \cite{Modular}]
\label{unimod-teo}
The Van den Bergh dualizing module
\begin{equation}
\label{V-Ah}
V_{\hA} = HH^d (\hA, \hA \otimes \hA )
\end{equation}
of $\hA = (A[[\h]], *)$ is isomorphic to $\hA$ as a bimodule
if and only if the formal Poisson structure
$\pi$ (\ref{pi}) corresponding to the star-product
$*$ is unimodular.
\end{teo}
The proof of this criterion is based on the algebraic
geometry version of Theorem \ref{chains}.
This criterion was used in recent paper \cite{EG} by P. Etingof
and V. Ginzburg to show that a certain family of Calabi-Yau
algebras associated to del Pezzo surfaces has a non-trivial
center.

Another application of formality theorems to quantization of
unimodular
Poisson structures is given in \cite{CF} by A.S. Cattaneo and
G. Felder. In this paper they showed that if
$\pi$ is a unimodular Poisson structure on a smooth real
manifold $M$ then the deformation quantization algebra
$(C^{\infty}(M)[[\h]], *)$ corresponding to $\pi$ admits
the following trace functional
\begin{equation}
f \mapsto \int_M \,f\, \Vol_{\h}\,,
\end{equation}
where $f$ is a compactly supported function on $M$ and
$\Vol_{\h}$ is a formal series of top degree forms starting
with the volume form $\Vol$ satisfying the equation
$$
l_{\pi} \Vol = 0\,.
$$

Algebraic index theorems \cite{Pindex}, \cite{TT} for
general Poisson manifolds give another tool for extracting information
about Kontsevich's class of a star-product.
These theorems express the isomorphism between Hochschild
(resp. periodic cyclic) homology of a deformation
quantization algebra and the Poisson homology (resp. de Rham
cohomology) of the manifold or variety in terms
of characteristic classes.  For the lack of
space we do not give more details about these theorems here
and instead refer the reader to \cite{Pindex} and \cite{TT}.

Many interesting examples of Poisson manifolds
are obtained via reduction \cite{Cattaneo}, \cite{Stasheff}.
In papers \cite{Rel-FT},
\cite{Tomsk} a ``super''-version
of Kontsevich's formality theorem is considered
with the application to deformation
quantization of reduced spaces
in a fairly general situation.
In both papers \cite{Rel-FT},
\cite{Tomsk} it was noticed that in general there may be
obstructions to the construction of the star-product
on the reduced space. This is not surprising because, in general,
reduction procedure gives a geometric object which is not even a manifold.
The presence of this obstruction does not mean that
certain reduced spaces should be discarded.
It is rather an indication that our
formulation of the quantization problem for such ``spaces''
should be modified.

The question of functoriality in deformation quantization
is closely related to the above question on the reduction.
We suspect that the ideas from paper \cite{Borisov} may
shed some light on this question.

There are also very interesting applications of
formality theorems for Hochschild complexes to Lie theory.
It is Kontsevich's formality theorem \cite{K} which
helped to solve \cite{AM}, \cite{Sahi}, \cite{Sahi1}, \cite{Sahi11},
\cite{Torossian}
the long standing Kashiwara-Vergne conjecture \cite{KV}.
Furthermore, using Kontsevich's formality theorem in \cite{CT}
A.S. Cattaneo and C. Torossian generalized some of Lichnerowicz's
results for the commutativity of the algebra of
invariant differential operators and solved
a long standing problem posed by M. Duflo for the expression of invariant
differential operators on any symmetric spaces in exponential coordinates.
They also developed a new method to construct characters for algebras of invariant
differential operators.

\section{An example of a non-formal DG Lie algebra}
\label{ex-non-formal}
Let $A= \bbK[x^1, \dots, x^d]$ be a polynomial algebra
in $d$ variables over $\bbK$\,.
It is obvious that the DG Lie algebra structure on
$\sCbu(A)$ restricts to the truncated
Hochschild complex $C^{\ge 1}(A)$ of $A$\,.
Furthermore, due to commutativity of $A$ the Hochschild
differential vanishes on the degree zero cochains.
Therefore, the cohomology of the truncated
Hochschild complex $C^{\ge 1}(A)$ is the
vector space
$$
\wedge^{\ge 1} \Der(A)
$$
of polynomial polyvector fields on $\bbK^d$
of degrees $\ge 1$\,.

In this subsection we prove that
\begin{teo}
\label{non-formal}
If $d$ is even then the DG Lie algebra $C^{\ge 1}(A)$
is non-formal.
\end{teo}
{\bf Remark 1.} This theorem answers a question of the referee of
our paper \cite{BLT}.

~\\
{\bf Remark 2.} We suspect that the case of odd dimension can
be considered similarly with a help of a regular constant
Poisson structure of maximal rank.

~\\
{\bf Proof} goes by contradiction.
The formality of the DG Lie algebra  $C^{\ge 1}(A)$
would still imply a bijection between the equivalence classes of
star-products and the equivalence classes of formal Poisson structures
on the affine space $\bbK^d$\,.

Let
\begin{equation}
\label{non-F}
F:\wedge^{\ge 1} \Der(A) \leadsto C^{\ge 1}(A)
\end{equation}
be an $\Linf$ quasi-isomorphism from the
graded Lie algebra $\wedge^{\ge 1} \Der(A)$ to the
DG Lie algebra $C^{\ge 1}(A)$\,.

As we mentioned in Section \ref{f-versus-nonf}
the structure map
\begin{equation}
\label{F-1}
F_1 : \wedge^{\ge 1} \Der(A) \to C^{\ge 1}(A)
\end{equation}
is a quasi-isomorphism of cochain complexes.

In general, this quasi-isomorphism may differ from the standard
Hochschild-Kostant-Rosenberg inclusion
\begin{equation}
\label{HKR}
I_{HKR} :\wedge^{\ge 1} \Der(A) \hookrightarrow C^{\ge 1}(A)\,.
\end{equation}
by a coboundary term.

If this coboundary term is non-zero then applying
Lemma 1 from \cite{erratum} we modify the
$\Linf$ quasi-isomorphism (\ref{non-F}) in such a way that
$F_1$ will coincide with (\ref{HKR}). Thus, we may assume, without
loss of generality, that
\begin{equation}
\label{F1-HKR}
F_1 = I_{HKR}\,.
\end{equation}

The correspondence
between the formal Poisson structures (\ref{pi}) and the star-products
is given by the assignment
$$
\pi \,\,\mapsto \,\, *
$$
\begin{equation}
\label{star1}
a * b = a b + \sum_{k=1}^{\infty} \frac{1}{k!}
F_k(\pi, \pi, \dots, \pi)(a,b) \,.
\end{equation}

Furthermore, $\pi$ allows us to twist the quasi-isomorphism
(\ref{non-F}) to the $\Linf$ morphism
\begin{equation}
\label{F-pi}
F^{\pi} : \wedge^{\ge 1} \Der(A)[[\h]]
\leadsto  C^{\ge 1}(\hA)\,,
\end{equation}
where $\hA$ is the algebra $A[[\h]]$ with the star-product
(\ref{star1}) and the DG Lie algebra $\wedge^{\ge 1} \Der(A)[[\h]]$
is considered with the Lichnerowicz differential $[\pi,\,]_{SN}$\,.

According to \cite{thesis} or \cite{Ezra-Lie}
the structure maps of the twisted $\Linf$ quasi-isomorphism
$F^{\pi}$ are given by the formula
\begin{equation}
\label{F-pi-n}
F^{\pi}_n(\ga_1, \dots, \ga_n) =
\sum_{k=1}^{\infty}\frac{1}{k!}
F_{n+k} (\pi, \pi, \dots, \pi, \ga_1, \ga_2, \dots, \ga_n)\,.
\end{equation}

An obvious analog of Proposition \ref{twist}
allows us to conclude that the structure map of the first level
\begin{equation}
\label{q-iso-1}
F^{\pi}_1 :  (\wedge^{\ge 1} \Der(A)[[\h]], [\pi,\,]_{SN})
\stackrel{\sim}{\to}  C^{\ge 1}(\hA)
\end{equation}
is a quasi-isomorphism of cochain complexes.

Thus localizing (\ref{F-pi}) in $\h$ we get the
following $\Linf$ quasi-isomorphism of DG Lie algebras
\begin{equation}
\label{F-pi-h-1}
F^{\pi} :  (\wedge^{\ge 1} \Der(A)((\h)), [\pi,\,]_{SN})
\leadsto  C^{\ge 1}(\hA[\h^{-1}])\,.
\end{equation}
The zeroth cohomology of the DG Lie algebra
$(\wedge^{\ge 1} \Der(A)((\h)), [\pi,\,]_{SN})$ is the
Lie algebra of Poisson vector fields and the zeroth cohomology of
the DG Lie algebra $C^{\ge 1}(\hA[\h^{-1}])$
is the Lie algebra of derivations of $\hA[\h^{-1}]$\,.

Therefore, the $\Linf$ quasi-isomorphism (\ref{F-pi-h-1})
gives us an isomorphism from the Lie algebra of Poisson vector
fields of $\pi$ to the Lie algebra $\Der(\hA[\h^{-1}])$
of derivations of $\hA[\h^{-1}]$\,.

Since $d$ is even we may choose $\pi= \h \te$
where $\te$ is a non-degenerate constant bivector.

Let us denote by $*_{W}$ the Moyal-Weyl star-product which
quantizes $\te$
\begin{equation}
\label{Weyl}
a\, *_W\, b =
\exp\left(\frac{\h}{2} \te^{ij}
\frac{\pa}{\pa x^i} \frac{\pa}{\pa y^j}\right) a(x)b(y)\, \Big|_{x=y} \,.
\end{equation}

Due to (\ref{F1-HKR}) the equivalence class of
$*_W$ corresponds to an equivalence class of a formal
Poisson structure $\tpi$ which starts with $\h \te$:
\begin{equation}
\label{tpi}
\tpi = \h \te + \h^2 \tpi_2 + \h^3 \tpi_3 + \dots.
\end{equation}

Since $\te$ is a non-degenerate bivector the formal
Poisson structure $\tpi$ is equivalent to the original
Poisson structure $\pi= \h \te$. This statement can be easily
deduced from the fact that, in the symplectic case, the Poisson
cohomology is isomorphic to the de Rham cohomology.

Therefore, the Moyal-Weyl star-product $*_W$ is equivalent
to $*$ (\ref{star1})\,.
Hence the Lie algebra $\Der(\hA^{W})$ of derivations of the
Weyl algebra
\begin{equation}
\label{A-Weyl}
\hA^{W} = (A((\h)), *_{W})
\end{equation}
is isomorphic
to the Lie algebra of Poisson vector fields of $\h\te$\,.

It is known that the Lie algebra $\Der(\hA^{W})$ of
derivations of the Weyl algebra $\hA^W$ is isomorphic to
\begin{equation}
\label{Lie-q}
\Big(\, \big( A/\bbK \big) ((\h)), [\, , \,]_{*_W}\, \Big)\,,
\end{equation}
where $[a,b]_{*_W} = a *_W  b -  b *_W a $\,.
The desired isomorphism is defined by assigning to an element
$a\in  \big( A/\bbK \big) ((\h))$ the
corresponding inner derivation:
$$
a ~~\mapsto~~ [a,\,]_{*_W}\,.
$$
Here we identify the quotient
$A/\bbK$ with the ideal of polynomials vanishing at the
origin.

Similarly, the Lie algebra of Poisson vector fields
of $\h \te$ is isomorphic to
\begin{equation}
\label{Lie-class}
\Big(\, \big( A/\bbK \big) ((\h)), \{\, , \,\}_{\h\te}\, \Big)\,,
\end{equation}
where $\{a,b\}_{\h \te} = \h\,  i_{\te}\, d a\, d b$\,.
The desired isomorphism is defined by assigning to an element
$a\in \big( A/\bbK \big) ((\h))$ the corresponding Hamiltonian
vector field:
$$
a ~~\mapsto~~ \{ a,\, \}_{\h \te}\,.
$$

Thus we conclude that an $\Linf$ quasi-isomorphism $F$
(\ref{non-F}) would give us an isomorphism $T$ from
the Lie algebra (\ref{Lie-class}) to the Lie algebra
(\ref{Lie-q}) over the field $\bbK((\h))$\,.
Furthermore, equations (\ref{F1-HKR}) and (\ref{F-pi-n})
imply that this isomorphism $T$ satisfies the following
property
\begin{equation}
\label{property}
T(a) = a ~~ {\rm mod}~~ \h\,, \qquad \forall~~ a\in A/\bbK\,.
\end{equation}

The Lie bracket in (\ref{Lie-q}) has the form
\begin{equation}
\label{Weyl-Lie}
[a, b]_{*_W} =  \h \te^{ij}
\pa_{x_i} a \pa_{x^j} b +
\h^3 V(a,b) ~~ {\rm mod} ~~ \h^5 \,,
\end{equation}
where
\begin{equation}
\label{Vey}
V(a,b) = \frac{1}{24}\, \te^{i_1 j_1}_1
\te^{i_2 j_2}_1 \te^{i_3 j_3}_1 \,
\pa_{x^{i_1}}\pa_{x^{i_2}}\pa_{x^{i_3}}\,a \,\,
\pa_{x^{j_1}}\pa_{x^{j_2}}\pa_{x^{j_3}}\,b\,,
\end{equation}
$$
a,b \in \big( A/\bbK \big) ((\h))\,.
$$
From deformation theory it follows that map
$$
V: A/\bbK \otimes A/\bbK \to A/\bbK
$$
given by the formula (\ref{Vey}) is a cocycle
for the Lie algebra  $A/\bbK$ with the bracket
\begin{equation}
\label{bracket}
\{a,b\} = \te^{ij}\, \pa_{x^i} a \,  \pa_{x^j} b\,.
\end{equation}
The existence of the isomorphism $T$ satisfying
the property (\ref{property}) would imply that this cocycle
is trivial. In other words, there should exist
a linear map
$$
P:  A/\bbK \to A/\bbK
$$
such that
\begin{equation}
\label{trivial}
V(a,b) = P(\{a,b\}) - \{P(a), b\} -\{ a, P(b)\}\,,
\end{equation}
$$
a,b\in A/\bbK\,.
$$

It is not hard to see that if the cocycle $V$ (\ref{Vey})
is trivial in the general case then it is trivial in
the two-dimensional case with the canonical Poisson
bracket
\begin{equation}
\label{bracket-dv}
\{a,b\}(x,y) = \pa_{x} a(x,y) \pa_{y} b(x,y) -
\pa_{y} a(x,y) \pa_{x} b(x,y)\,,
\end{equation}
$$
a,b \in \bbK[x,y]/\bbK \,.
$$
Thus we may restrict ourselves to the two-dimensional case with
the canonical Poisson bracket (\ref{bracket-dv}).

In this case the cocycle $V$ reads
\begin{equation}
\label{Vey-dv}
V(a,b) = \frac{1}{24}\big(
\pa^3_x (a) \pa^3_y b
-3 \pa^2_x \pa_y (a) \pa_x \pa^2_y (b)
+ 3 \pa_x \pa^2_y (a) \pa^2_x \pa_y (b)
- \pa^3_y (a) \pa^3_x (b)
\big)\,,
\end{equation}
$$
a,b \in \bbK[x,y]/\bbK \,.
$$

For $a=x$ and $b=y$ equation (\ref{trivial}) implies that
$$
\pa_x P(x) + \pa_y P(y) \in \bbK\,.
$$

Let $c_0 = \displaystyle \frac{1}{2} (\pa_x P(x) + \pa_y P(y))$\,.
Then, setting $Q_1 = P(x)- c_0 x$ and $Q_2 = P(y)- c_0 y$\,,
we get
$$
\pa_x Q_1 + \pa_y Q_2 =0\,.
$$
Hence, there exists a polynomial $Q\in \bbK[x,y]$ such that
$Q_{1}= \{Q, x\}$ and $Q_2=\{Q , y\}$\,.

Thus, adjusting $P$ by a Hamiltonian vector field, we
reduce it to the form in which
\begin{equation}
\label{P-linear}
P(x) = c_0 x\,, \qquad P(y) = c_0 y\,.
\end{equation}

Substituting to equation (\ref{trivial}) quadratic
monomials for $a$ and linear for $b$ we deduce that
$$
\begin{array}{c}
P(x^2) = c_{11} x + c_{12} y \,, \\[0.3cm]
P(xy) = c_{21} x + c_{22} y \,, \\[0.3cm]
P(y^2) = c_{31} x + c_{32} y \,,
\end{array}
$$
where $c_{ij}\in \bbK$\,.
Next substituting quadratic monomials for $a$ and $b$
we deduce that $c_{12} = c_{31} = 0 $\,, $c_{11}= 2 c_{22}$
and $c_{32}= 2 c_{21}$\,. Thus
\begin{equation}
\label{P-quad}
\begin{array}{c}
P(x^2) = 2 c_{y} x\,,  \\[0.3cm]
P(xy) = c_{x} x + c_{y} y\,, \\[0.3cm]
P(y^2) = 2 c_{x} y
\end{array}
\end{equation}
for some constants $c_x, c_y\in \bbK$\,.

Adjusting $P$ by the Hamiltonian vector field
$$
\{c_x x - c_y y, \,\}
$$
we kill the right hand sides in (\ref{P-quad})\,.

Thus we may assume that
\begin{equation}
\label{P-linear-quad}
\begin{array}{c}
\begin{array}{cc}
P(x) = c_0 x\,, & P(y) = c_0 y\,,
\end{array} \\[0.3cm]
P(x^2) = P(xy) = P(y^2) = 0\,.
\end{array}
\end{equation}

Next, plugging in cubic monomials $x^3$, $x^2 y$, $x y^2$, $y^3$
for $a$ and linear monomials for $b$ in (\ref{trivial}),
we deduce that
\begin{equation}
\label{P-cubic}
\begin{array}{c}
P(x^3) = b_{11}x + b_{12}y - c_0 x^3 \,,  \\[0.3cm]
P(x^2 y) = b_{21}x + b_{22}y - c_0 x^2 y \,, \\[0.3cm]
P(x y^2) = b_{31}x + b_{32}y - c_0 x y^2 \,, \\[0.3cm]
P(y^3) = b_{41}x + b_{42}y - c_0 y^3\,,
\end{array}
\end{equation}
where $b_{ij} \in \bbK$\,.

To get a further restriction we substitute cubic monomials
for $a$ and quadratic for $b$ in (\ref{trivial}). We get that
all the coefficients $b_{ij}$ should vanish. Thus
\begin{equation}
\label{P-cubic1}
\begin{array}{cc}
P(x^3) = - c_0 x^3 \,, & P(x^2 y) = - c_0 x^2 y \,,\\[0.3cm]
P(x y^2) =  - c_0 x y^2 \,, & P(y^3) = - c_0 y^3\,,
\end{array}
\end{equation}

Every quartic monomial can be written as a Poisson
bracket of two cubic monomials. Using this observation
and equation (\ref{trivial}) for $a$ and $b$ being
cubic monomials we deduce that
$$
P(x^n y^k) = - 2 c_0 x^n y^k
$$
whenever $n+k=4$\,.

To get the desired contradiction we, first, set $a= x^4$
and $b=y^3$ in (\ref{trivial}) and get
$$
6 x - (12 P(x^3 y^2) + 2 c_0\{x^4, y^3\} + c_0 \{x^4 , y^3 \}) \in \bbK
$$
or equivalently
\begin{equation}
\label{raz}
12 P(x^3 y^2) - 6 x + 3 c_0 x^3 y^2 \in \bbK \,.
\end{equation}

Second, plugging $a= x^3 y$
and $b= x y^2$ into (\ref{trivial}) we get
$$
-\frac{3}{2} x - (5 P(x^3 y^2) + 2 c_0\{x^3 y, x y^2\} + c_0 \{x^3 y, x y^2 \})
$$
or equivalently
\begin{equation}
\label{dva}
5 P(x^3 y^2) - \frac{3}{2} x + 3 c_0 x^3 y^2 \in \bbK\,.
\end{equation}
The inclusion in (\ref{raz}) clearly contradicts to
the inclusion in (\ref{dva}) and the theorem follows. $\Box$

~\\
{\bf Remark 1.} The same expression (\ref{Vey}) defines
a cocycle for the Lie algebra $A$ with the bracket
$\{\,,\,\}_{\te}$\,. It was shown by J. Vey in \cite{Vey} that this
cocycle is non-trivial. Here we had to work with the quotient
$A/\bbK$ and this is why we had to redo the computation of J. Vey
taking into account this modification.

~\\
{\bf Remark 2.} It makes sense to consider a modification
of the Weyl algebra which is defined over $\bbK$\,.
This is the algebra $A^W$  generated
by $x^i$'s satisfying the relations
$$
x^i \, x^j - x^j \, x^i = \te^{ij}\,,
$$
where $\te^{ij}$ is as above a non-degenerate antisymmetric
constant matrix.  It is known \cite{Belov-M} that the Lie algebra of
derivation of $A^W$ is not isomorphic
to the Lie algebra of the derivation of the corresponding
Poisson algebra $\bbK[x^1, \dots, x^d]$ with the bracket
$\{\, ,\, \}_{\te}$\,. This fact is mentioned in \cite{Belov-M}
as a negative evidence for the conjecture about the automorphisms
of the Weyl algebra.

~\\

\noindent\textsc{Department of Mathematics,
University of California at Riverside, \\
900 Big Springs Drive,\\
Riverside, CA 92521, USA \\
\emph{E-mail address:} {\bf vald@math.ucr.edu}}

~\\

\noindent\textsc{Mathematics Department,
Northwestern University, \\
2033 Sheridan Rd.,\\
Evanston, IL 60208, USA \\
\emph{E-mail addresses:} {\bf tamarkin@math.northwestern.edu},
{\bf tsygan@math.northwestern.edu}}


\begin{thebibliography}{99}

\bibitem{AM} A. Alekseev and E. Meinrenken, On the Kashiwara-Vergne
conjecture, Invent. Math. {\bf 164} (2006) 615--634.


\bibitem{Sahi} M. Andler, A. Dvorsky, and S. Sahi, Deformation
quantization and invariant distributions,
{\it C. R. Acad. Sci. Paris S\'er. I Math.} {\bf 330}, 2 (2000) 115--120.

\bibitem{Sahi1} M. Andler, A. Dvorsky, and S. Sahi, Kontsevich
quantization and invariant distributions on Lie groups,
{\it Ann. Sci. Ecole Norm. Sup.} (4) {\bf 35}, 3 (2002) 371--390.

\bibitem{Sahi11} M. Andler, S. Sahi, and C. Torossian,
Convolution of invariant distributions: proof of the Kashiwara-Vergne conjecture,
{\it Lett. Math. Phys.} {\bf 69} (2004) 177--203.


\bibitem{signs} D. Arnal, D. Manchon, and M. Masmoudi,
Choix des signes pour la formalit\'e de M. Kontsevich,
Pacific J. Math. {\bf 203}, 1 (2002) 23--66.

\bibitem{Babenko} I. K. Babenko and I. A. Taimanov,
Massey products in symplectic manifolds, Sb. Math. {\bf 191},
7-8 (2000) 1107--1146.

\bibitem{Bayen}  F. Bayen, M. Flato, C. Fronsdal, A. Lichnerowicz, and D.
Sternheimer,  Deformation theory and quantization. I. Deformations
of symplectic structures, Ann. Phys. (N.Y.), {\bf 111} (1978)
61;\\ Deformation theory and quantization, II. Physical
applications, Ann. Phys. (N.Y.), {\bf 110} (1978) 111.


\bibitem{Belov-M} A. Belov-Kanel and M. Kontsevich,
Automorphisms of the Weyl algebra,
Lett. Math. Phys., {\bf 74}, 2 (2005) 181--199;
arXiv:math/0512169.


\bibitem{Ber}  F.A. Berezin, Quantization, Izv. Akad. Nauk., {\bf 38}
(1974) 1116-1175;\\ General concept of quantization, Commun. Math.
Phys., {\bf 40} (1975) 153-174.


\bibitem{BF} C. Berger and B. Fresse,
Combinatorial operad actions on cochains,
Math. Proc. Cambridge Philos. Soc.,
{\bf 137}, 1  (2004) 135--174.

\bibitem{BM} C. Berger and I. Moerdijk,
Axiomatic homotopy theory for operads,
Comment. Math. Helv. {\bf 78}, 4 (2003) 805--831;
arXiv:math/0206094

\bibitem{Board-V} J.M. Boardmann and R.M. Vogt,
Homotopy invariant algebraic structures on topological spaces,
Springer-Verlag, Berlin, 1973, Lect. Notes in Math., Vol. 347.

\bibitem{Borisov} D.V. Borisov, $G_{\infty}$-structure on the
deformation complex of a morphism, J. Pure Appl. Algebra
{\bf 210}, 3 (2007) 751--770.

\bibitem{JLB} J.-L. Brylinski,
A differential complex for Poisson manifolds,
J. Diff. Geom. {\bf 28}, 1  (1988) 93--114.

\bibitem{BZ} J.-L. Brylinski and G. Zuckerman,
The outer derivation of a complex Poisson manifold,
J. Reine Angew. Math. {\bf 506} (1999) 181--189.

\bibitem{CH} D. Calaque and G. Halbout, Weak quantization of
Poisson structures, arXiv:0707.1978.


\bibitem{CV} D. Calaque and M. Van den Bergh,
Hochschild cohomology and Atiyah classes,
arXiv:0708.2725.

\bibitem{CV-Ginfty} D. Calaque and M. Van den Bergh,
Global formality at the $G_\infty$-level,
arXiv:0710.4510.

\bibitem{CR-cap} D. Calaque and C. A. Rossi,
Shoikhet's Conjecture and Duflo Isomorphism on (Co)Invariants,
SIGMA {\bf 4} (2008) 060; arXiv:0805.2409.


\bibitem{Cald} A. C\u{a}ld\u{a}raru, The Mukai pairing. II.
The Hochschild-Kostant-Rosenberg isomorphism,
Adv. Math.,  {\bf 194}, 1  (2005) 34--66.

\bibitem{Cattaneo} A.S. Cattaneo, Deformation quantization and reduction.
Poisson geometry in mathematics and physics, 79--101, Contemp. Math., 450,
Amer. Math. Soc., Providence, RI, 2008.


\bibitem{CF-path} A.S. Cattaneo and G. Felder,
A path integral approach to the Kontsevich quantization formula,
Commun. Math. Phys. {\bf 212}, 3 (2000) 591--611;
arXiv:math/9902090.


\bibitem{CFT} A.S. Cattaneo, G. Felder, and L. Tomassini,
From local to global deformation quantization of Poisson
manifolds, Duke Math. J., {\bf 115}, 2  (2002) 329-352;
math.QA/0012228.


\bibitem{Rel-FT} A.S. Cattaneo and G. Felder,
Relative formality theorem and quantisation of coisotropic
submanifolds,  Adv. Math. {\bf 208}, 2 (2007) 521--548.

\bibitem{CF} A. S. Cattaneo and G. Felder,
Effective Batalin-Vilkovisky theories, equivariant configuration spaces
and cyclic chains, arXiv:0802.1706.

\bibitem{CT} A. S. Cattaneo and C. Torossian, Quantification pour les paires
symetriques et diagrammes de Kontsevich, arXiv:math/0609693.

\bibitem{Boris-book}
J. Cuntz, G. Skandalis, and B. Tsygan.
{\it Cyclic homology in non-commutative geometry.}
Encyclopaedia of Mathematical Sciences, 121.
Operator Algebras and Non-commutative Geometry,
II. Springer-Verlag, Berlin, 2004.


\bibitem{DGT} Yu. Daletski, I. Gelfand, and
B. Tsygan, On a variant of noncommutative
geometry, Soviet Math. Dokl. {\bf 40}, 2
(1990) 422--426.

\bibitem{DGMS} P. Deligne, P. Griffiths, J. Morgan, and D. Sullivan,
Real homotopy theory of K\"ahler manifolds,
Invent. Math. {\bf 29}, 3 (1975) 245--274.



\bibitem{CEFT} V.A. Dolgushev,
Covariant and Equivariant Formality Theorems,
Adv. Math., {\bf 191}, 1 (2005) 147--177;
arXiv:math/0307212.

\bibitem{FTHC} V.A. Dolgushev,
A Formality Theorem for Hochschild Chains,
Adv. Math.  {\bf 200}, 1  (2006) 51--101;
math.QA/0402248.

\bibitem{thesis} V.A. Dolgushev,
A Proof of Tsygan's formality conjecture
for an arbitrary smooth manifold, PhD thesis,
MIT; math.QA/0504420.

\bibitem{Modular} V. A. Dolgushev, The Van den Bergh duality and
the modular symmetry of a Poisson variety, accepted to
Selecta Math.; arXiv:math/0612288.

\bibitem{erratum} V.A. Dolgushev,
Erratum to: "A Proof of Tsygan's Formality Conjecture for an
Arbitrary Smooth Manifold", arXiv:math/0703113.


\bibitem{Pindex} V.A. Dolgushev and V.N. Rubtsov,
An algebraic index theorem for Poisson manifolds,
arXiv:0711.0184.


\bibitem{BLT} V. Dolgushev, D. Tamarkin and
B. Tsygan, The homotopy Gerstenhaber algebra of Hochschild
cochains of a regular algebra is formal, J. Noncommut. Geom.
{\bf 1}, 1 (2007) 1--25;
arXiv:math/0605141.


\bibitem{Cepochki}  V.A. Dolgushev, D.E. Tamarkin, and
B.L. Tsygan, Formality of the homotopy calculus algebra of
Hochschild (co)chains, arXiv:0807.5117.



\bibitem{Drinfeld} V.G. Drinfeld, Quasi-Hopf algebras,
Leningrad Math. J. {\bf 1}, 6 (1990) 1419--1457.


\bibitem{EG} P. Etingof and V. Ginzburg,
Noncommutative del Pezzo surfaces and Calabi-Yau algebras,
arXiv:0709.3593.



\bibitem{Fedosov} B.V. Fedosov, A simple geometrical construction
of deformation quantization, J. Diff. Geom. {\bf 40} (1994)
213--238.

\bibitem{GF} I.M. Gelfand and D.V. Fuchs,
Cohomology of the algebra of formal vector
fields, Izv. Akad. Nauk., Math. Ser.
{\bf 34} (1970) 322--337 (In Russian).

\bibitem{G-Kazh} I.M. Gelfand and D.A. Kazhdan,
Some problems of differential geometry and
the calculation of cohomologies of Lie algebras of
vector fields, Soviet Math. Dokl., {\bf 12},
5 (1971) 1367-1370.


\bibitem{G} M. Gerstenhaber, The cohomology structure of
an associative ring,  Annals of Math., {\bf 78} (1963)
267--288.

\bibitem{Ezra-higher} E. Getzler, A Darboux theorem for Hamiltonian operators
in the formal calculus of variations, Duke Math. J. {\bf 111}, 3 (2002)
535--560.

\bibitem{Ezra-Lie} E. Getzler,
Lie theory for nilpotent L-infinity algebras,
to appear in Ann. Math.; arXiv:math/0404003

\bibitem{GJ} E. Getzler and J.D.S. Jones,
Operads, homotopy algebra and iterated integrals
for double loop spaces, hep-th/9403055.


\bibitem{Vitya} V. Ginzburg, Calabi-Yau algebras,
math.AG/0612139.


\bibitem{GK} V. Ginzburg and M. Kapranov,
Koszul duality for operads,  Duke Math. J.
{\bf 76}, 1 (1994) 203--272.


\bibitem{Goldman-M} W. Goldman and J. Millson. The deformation
theory of representation of fundamental groups in compact
K\"ahler manifolds. {\it Publ. Math. I.H.E.S.},
{\bf 67} (1988) 43-96.

\bibitem{Halperin-Stasheff} S. Halperin and J. Stasheff,
Obstructions to homotopy equivalences, Adv. Math. {\bf 32}, 3 (1979)
233--279.

\bibitem{Henriques} A. Henriques, Integrating $L_\infty$-algebras,
Compos. Math. {\bf 144}, 4 (2008) 1017--1045.

\bibitem{Hinich-op} V. Hinich, Homological algebra of homotopy
algebras, Comm. Algebra {\bf 25}, 10 (1997) 3291--3323.

\bibitem{Hinich} V. Hinich, Tamarkin's proof of
Kontsevich formality theorem, Forum Math. {\bf 15},
4 (2003) 591--614; math.QA/0003052.


\bibitem{HS} V. Hinich and V. Schechtman, Homotopy Lie algebras,
I.M. Gelfand Seminar, Adv. Sov. Math., {\bf 16}, 2 (1993) 1-28.


\bibitem{Ikeda} N. Ikeda, Two-dimensional gravity and nonlinear gauge theory,
Ann. Phys. {\bf 235} (1994) 435–-464.


\bibitem{KV} M. Kashiwara and M. Vergne, The Campbell-Hausdorff
formula and invariant hyperfunctions, Invent. Math.
{\bf 47} (1978) 249–-272.

\bibitem{K} M. Kontsevich, Deformation quantization of
Poisson manifolds, Lett. Math. Phys.,
{\bf 66} (2003) 157-216; q-alg/9709040.

\bibitem{K1} M. Kontsevich, Formality Conjecture,
{\it D. Sternheimer et al. (eds.), Deformation Theory and Symplectic
Geometry}, Kluwer 1997, 139 -- 156.

\bibitem{K-alg} M. Kontsevich, Deformation quantization of algebraic varieties.
Mosh\'e Flato memorial conference 2000, Part III (Dijon).
Lett. Math. Phys. {\bf 56}, 3 (2001) 271--294.

\bibitem{K-oper} M. Kontsevich,
Operads and motives in deformation quantization,
Lett. Math. Phys., {\bf 48} (1999) 35--72.


\bibitem{K-Soi} M. Kontsevich and Y. Soibelman,
Deformations of algebras over operads and the Deligne conjecture,
{\it Proceedings of the Mosh\'e Flato Conference
Math. Phys. Stud.} {\bf 21},
255--307, Kluwer Acad. Publ., Dordrecht, 2000.

\bibitem{K-Soi1} M. Kontsevich and Y. Soibelman,
Notes on A-infinity algebras, A-infinity categories
and non-commutative geometry. I, math.RA/0606241.


\bibitem{Koszul} J.L. Koszul,
Crochet de Schouten-Nijenhuis et cohomologie,
Ast\'erisque  (1985)  Numero Hors Serie, 257--271.


\bibitem{LS} T. Lada and J. Stasheff, Introduction to SH Lie
algebras for physicists, Intern. J. Theor Phys. {\bf 32}, 7 (1993)
1087-1103.

\bibitem{LV} P. Lambrechts and I. Volic,
Formality of the little N-disks operad, arXiv:0808.0457.


\bibitem{Lich} A. Lichnerowicz,
Les vari\'et\'es de Poisson et leurs alg\`ebres de Lie
associ\'ees, J. Diff. Geom., {\bf 12}, 2  (1977)
253--300.

\bibitem{Loday} J.- L. Loday, {\it Cyclic Homology},
Grundlehren der mathematischen Wissenschaften, 301.
Springer-Verlag, Berlin, 1992.


\bibitem{Tomsk} S.L. Lyakhovich and A.A. Sharapov,
BRST theory without Hamiltonian and Lagrangian,
J. High Energy Phys. {\bf 3} (2005) 011, 22 pp.

\bibitem{Torossian-cup} D. Manchon and C. Torossian,
Cohomologie tangente et cup-produit pour la quantification de Kontsevich,
Ann. Math. Blaise Pascal {\bf 10}, 1 (2003) 75--106.

\bibitem{Markl} M. Markl, Models for operads, Comm. Algebra {\bf 24} (1996)
1471--1500.

\bibitem{May} J.P. May, Infinite loop space theory,
Bull. Amer. Math. Soc.,  {\bf 83}, 4  (1977) 456--494.


\bibitem{May-Massey} J.P. May, Matrix Massey products,
J. Algebra {\bf 12} (1969) 533--568.


\bibitem{M-Smith} J. E. McClure and J. H. Smith,
A solution of Deligne's Hochschild cohomology conjecture,
Recent progress in homotopy theory (Baltimore, MD, 2000),
{\it Amer. Math. Soc.,  Contemp. Math.}, {\bf 293}, 153--193;
math.QA/9910126.


\bibitem{Q} D. Quillen, Rational homotopy theory, Annals
of Math., {\bf 90}, 2 (1969) 205--295.


\bibitem{Thomas} P. Schaller and T. Strobl,
Poisson structure induced (topological) field theories,
Modern Phys. Lett. {\bf A9}, 33 (1994) 3129--3136.

\bibitem{SS-MC} M. Schlessinger and J. Stasheff,
Deformation theory and rational homotopy type, University of
North Carolina preprint, 1979.

\bibitem{Massey-Lie} M. Schlessinger and J. Stasheff,
The Lie algebra structure of tangent cohomology and deformation theory,
J. Pure Appl. Algebra {\bf 38}, 2-3 (1985) 313--322.

\bibitem{Sh} B. Shoikhet,
A proof of the Tsygan formality conjecture for chains,
Adv. Math., {\bf 179}, 1 (2003) 7--37;
math.QA/0010321.

\bibitem{Stasheff} J. Stasheff, Homological reduction of constrained Poisson
algebras, J. Diff. Geom. {\bf 45} (1997) 221–-240.

\bibitem{Swan} R.G. Swan, Hochschild cohomology of
quasiprojective schemes, J. Pure Appl. Algebra {\bf 110}, 1
(1996) 57--80.

\bibitem{Dima-Proof} D. Tamarkin,
Another proof of M. Kontsevich formality theorem,
math.QA/9803025.

\bibitem{Dima-Disk} D. Tamarkin,
Formality of chain operad of little discs,
Lett. Math. Phys.  {\bf 66}, 1-2  (2003) 65--72;
math.QA/9809164.

\bibitem{Dima-DG} D. Tamarkin, What do DG categories
form?  Compos. Math. {\bf 143}, 5 (2007) 1335--1358;
math.CT/0606553.

\bibitem{TT} D. Tamarkin and B. Tsygan,
Cyclic formality and index theorems,
Talk given at the Mosh\'e Flato Conference (2000),
Lett. Math. Phys. {\bf 56}, 2  (2001) 85--97.

\bibitem{Torossian} C. Torossian, Sur la conjecture combinatoire de
Kashiwara-Vergne, J. Lie Theory {\bf 12}, 2 (2002) 597--616.

\bibitem{Tsygan} B. Tsygan,
Formality conjectures for chains,
Differential topology, infinite-dimensional Lie algebras,
and applications. 261--274,
{\it Amer. Math. Soc. Transl.} Ser. 2, 194,
Amer. Math. Soc., Providence, RI, 1999.


\bibitem{VB} M. Van Den Bergh, A Relation between Hochschild
Homology and Cohomology for Gorenstein Rings, Proc. Amer. Math.
Soc. {\bf 126}, 5 (1998) 1345-1348;\\
Erratum to  ``A Relation between Hochschild
Homology and Cohomology for Gorenstein Rings'',
Proc. Amer. Math. Soc. {\bf 130}, 9 (2002) 2809-2810.


\bibitem{VdB-glob} M. Van den Bergh, On global deformation quantization
in the algebraic case, J. Algebra {\bf 315}, 1 (2007) 326--395.


\bibitem{Vey} J. Vey. D\'eformation du crochet de
Poisson sur une vari\'et\'e
symplectique, Comment. Math. Helv.,
{\bf 50} (1975) 421--454.


\bibitem{Sasha-V} A. Voronov, Quantizing Poisson manifolds,
{\it Perspectives on quantization} (South Hadley, MA, 1996), 189--195,
Contemp. Math., {\bf 214}, Amer. Math. Soc., Providence, RI, 1998.


\bibitem{Sasha-Proof} A.A. Voronov, Homotopy Gerstenhaber
algebras, {\it Proceedings of the Mosh\'e Flato Conference
Math. Phys. Stud.}, {\bf 22}, 307-331.
Kluwer Acad. Publ., Dordrecht, 2000.


\bibitem{Alan} A. Weinstein, The modular automorphism group
of a Poisson manifold,
J. Geom. Phys.  {\bf 23}, 3-4  (1997) 379--394.


\bibitem{W} T. Willwacher, Formality of cyclic chains,
arXiv:0804.3887.


\bibitem{Y} A. Yekutieli, The continuous Hochschild cochain
complex of a scheme, Canad. J. Math. {\bf 54}, 6 (2002) 1319--1337.

\bibitem{Ye} A. Yekutieli, Deformation quantization in
algebraic geometry, Adv. Math. {\bf 198}, 1 (2005) 383--432;
math.AG/0310399.


\end{thebibliography}
\end{document}